\documentclass{article}
\usepackage{amsmath}
\usepackage{amsthm}
\usepackage{amssymb}
\usepackage{amsfonts}
\usepackage{mathrsfs}
\usepackage{enumerate}
\usepackage{authblk}
\usepackage[normalem]{ulem} 
\usepackage[dvipsnames]{xcolor}
\usepackage[utf8]{inputenc}
\usepackage[english]{babel}
\usepackage[ top = 3 cm, bottom = 2 cm]{geometry} 
\usepackage{graphicx}

\usepackage[colorlinks,citecolor=red,pagebackref,hypertexnames=false]{hyperref}
\usepackage{pdfsync}

\newcommand{\ep}{\varepsilon}
\newcommand{\n}[1]{\mathscr{#1}}
\newcommand{\m}[1]{\mathcal{#1}}
\newcommand{\bb}[1]{\mathbb{#1}}
\newcommand{\ra}{\rightarrow}

\newcommand{\D}{\Delta}
\DeclareMathOperator*{\esssup}{esssup}

\numberwithin{equation}{section}
\theoremstyle{plain}
\newtheorem{theorem}{Theorem}
\newtheorem{lemma}[theorem]{Lemma}
\newtheorem{proposition}[theorem]{Proposition}
\newtheorem{corollary}[theorem]{Corollary}

\theoremstyle{definition}
\newtheorem{definition}[theorem]{Definition}

\begin{document}

\title{Optimal Control of Singular Parabolic PDEs Modeling Multiphase Stefan-type Free Boundary Problems}
\author{Ugur G. Abdulla\thanks{abdulla@fit.edu} \ } 
\author{Evan Cosgrove\thanks{ecosgrove2011@my.fit.edu}}
\affil{Department of Mathematical Sciences, Florida Institute of Technology,\\ Melbourne, Florida 32901}

\maketitle

\abstract{Optimal control of the singular nonlinear parabolic PDE which is a distributional formulation of multidimensional and multiphase Stefan-type
free boundary problem is analyzed. Approximating sequence of finite-dimensional 
optimal control problems is introduced via finite differences. Existence of the optimal control and the convergence of the sequence of discrete optimal
control problems both with respect to functional and control is proved. In particular, convergence of the method of finite differences, and existence, uniqueness and stability estimations are established for the singular PDE problem under minimal regularity assumptions on the coefficients.
}\\

{{\bf Key words}: optimal control, singular PDE, multiphase free boundary problem, nonlinear parabolic PDE with discontinuous coefficient, discrete optimal control problem, convergence in functional, convergence in control.}

{{\bf AMS subject classifications}: 35R30, 35R35, 35K20, 35Q93, 49J20, 65M06, 65M12} 

\newpage
\section{Introduction}

\subsection{Optimal Control Problem}

Let $d\in\bb N, \Omega\subset \bb R^d$ be a bounded domain with Lipschitz boundary, $T>0$, $D:=\Omega\times(0,T]$ and $v^1<v^2<\cdots<v^m$ are given real numbers. Consider singular PDE problem:
\begin{align}
    \frac{\partial \beta(v)}{\partial t}-\mathcal{L}v-f(x,t) \ni 0, &\qquad (x,t)\in D,\label{PDE}\\
    v(x,0)=\Phi(x) \label{vphifuture}, &\qquad x \in \Omega,\\
    v\vert_S = 0,&\qquad 0<t\leq T,  \label{vgfuture}
\end{align}
where $\beta(\cdot)$ is a maximal monotone graph of the form
\begin{equation}\label{beta}
    \beta(y)=
    \left\{
    \begin{array}{l}
    \beta_j(y)+\sum\limits_{i=0}^{j-1}\nu_i, \quad\text{for} \ v^{j-1}<y<v^j,\\
    $\Big[$\beta_j(v^j)+\sum\limits_{i=0}^{j-1}\nu_i,\beta_j(v^j)+\sum\limits_{i=0}^{j}\nu_i $\Big]$, \quad\text{for} \  \ y=v^j,\\
    \beta_{j+1}(y)+\sum\limits_{i=0}^{j}\nu_i, \quad\text{for} \ v^j<y<v^{j+1}; \ j=1,2,...,m
    \end{array}\right.
\end{equation}
with a given positive constants $\nu_j, j=1,...,m$; $\nu_0=0$, $v^0=-\infty$, $v^{m+1}=+\infty$;  $\beta_i(\cdot), i=1,...,m+1$ are monotone increasing Lipschitzian functions in their respective domain of definition, $\beta_j(v^j)=\beta_{j+1}(v^j), j=1,...,m$, 
    \begin{equation}\label{bbar} 
         \beta_j^{\prime}(y)\geq \bar{b}>0, \ j=1,...,m+1;
    \end{equation}
and $\mathcal{L}$ is an elliptic operator
\begin{gather}
	\mathcal{L}v = \sum\limits_{i=1}^{d}(a_{i}(x,t)v_{x_i}+b_i(x,t)v)_{x_i}-\sum\limits_{i=1}^d c_i(x,t)v_{x_i} - r(x,t)v
\end{gather}
with bounded and measurable coefficients $a_i,b_i,c_i,r$ and 
\begin{equation}\label{apositive}
     a_{i}(x,t) \geq a_0 >0, \quad i=\overline{1,d}, \quad  \text{a.e.} \ (x,t)\in D.
\end{equation}

Singular PDE problem \eqref{PDE}-\eqref{vgfuture} is a distributional formulation of the multiphase Stefan problem if elliptic operator $\mathcal{L}$ coincides with $\Delta$ \cite{Oleinik, Kamenomostskaya2, LSU}. In the physical context, $v(x,t)$ is a temperature distribution, $f(x,t)$ is a density of heat sources, $\Phi(x)$ is an initial temperature, $v^j$'s are phase transition temperatures; $\beta_j^{\prime}(v), v^j<v<v^{j+1}, j=0,1,...,m$ express heat conductivities in each phase, positive constants $\nu^j, j=1,...,m$ characterize latent heat of fusion during phase transition, and the coefficients $a_i,b_i,c_i,r$ characterize anisotropic properties of the media. The classical case $m=1, v^1=0, \mathcal{L}=\Delta$ is a two-phase Stefan problem describing melting of the ice or freezing of the water \cite{Friedman1,Meyrmanov}, whereas more sophisticated applications include  biomedical problem about the laser ablation of biomedical tissues, which motivates general elliptic operator $\mathcal{L}$. 

Consider optimal control problem on the minimization of the functional 
\begin{equation}\label{functional}
\n J(f)=\Vert v|_{\Omega\times\{t=T\}}-\Gamma\Vert^2_{L_2(\Omega)}
\end{equation}
on a control set
\[\n F^R =\left\{f\in L_{\infty}(D):~\Vert f\Vert_{L_{\infty}(D)}\leq R\right\}.\]
where $R>0$ and $\Gamma\in L_2(\Omega)$ are given, and $v=v(x,t; f)$ is a solution of the singular PDE problem \eqref{PDE}-\eqref{vgfuture}. Furthermore, this optimal control problem will be referred to as \emph{Problem $\m I$}. Motivation for the \emph{Problem $\m I$}  arises in many applications, such as a modeling and control of biomedical engineering problem about the laser ablation of biomedical tissues \cite{Abdulla7,Abdulla8}, preventing aerodynamic stall in aircrafts due to in-flight ice accretion \cite{Myers}, etc. The goal is to find optimal choice of the density of the sources which minimizes the mismatch of the temperature distribution at the final moment with the desired temperature profile.
 
In this paper existence of the optimal control, and convergence of the finite-difference approximations of the \emph{Problem $\m I$} both with respect to functional and control will be proved. In particular, convergence of the finite-difference approximations of the singular PDE problem, and existence, uniqueness and stability results will be established.

The transformation of the multiphase Stefan problem to singular PDE problem was first introduced in \cite{Oleinik}. Existence and uniqueness of the weak solution of the Stefan problem in a new formulation was established in \cite{Oleinik, Kamenomostskaya2,LSU}. In \cite{DiBenedetto1,DiBenedetto2} it is proved that weak solutions are H\"{o}lder continuous for general nonlinear elliptic operators $\mathcal{L}$. Weak solutions to the two phase Stefan problem were proven to be continous in \cite{Caffarelli}.

Besides of its own importance, optimal control framework is commonly used variational method for solving inverse Stefan-type free boundary problems. Motivated with different class of applications, historically, optimal control and inverse Stefan-type free boundary problems were developed in two different directions depending whether or not free boundary is known. 
One-dimensional and one phase inverse Stefan problem (ISP) was first considered in \cite{Cannon3}, where missing heat flux on the fixed boundary must be found with additional measurement of the phase transition boundary - the problem being reminiscent of the characteristic Cauchy problem for the heat equation. In \cite{BudakVasileva1, BudakVasileva2}, a variational approach was implemented to solve the ISP. The first formulation of the one-dimensional and one phase ISP with unknown free boundary as an optimal control problem appeared in \cite{Vasilev}, where existence of the optimal control was proved. In \cite{Yurii}, Fr\'{e}chet differentiability and convergence of finite difference schemes was proved, and Tikhonov regularization was implemented. Over the last half-century research on optimal control and ISP with given phase transition boundaries were addressed in \cite{Alifanov,Bell,Budak,Cannon,Carasso,Ewing1,Ewing2,Hoffman,sprekels1,Sherman,pawlow1,pawlow2,dunbar1,Goldman,Frankel,Zabaras1,Zabaras2,Zabaras3}, whereas the problems with unknown phase transition boundaries are addressed in \cite{Baumeister,Fasano,Hoffman1,Hoffman2,Jochum2,Jochum1,Knabner,Lurye,Nochetto,Niezgodka,Primicero,Sagues,Talenti,Goldman,Gol'dman,Hinze1,Hinze2,herzog}. 
Summarizing, the main methods to solve inverse Stefan-type free boundary problems were based on variational formulation in optimal control framework, method of quasisolutions and Tikhonov regularization to address ill-posedness in terms of errors of measurements, Frechet differentiability and gradient descent type iterative algorithms for the numerical solution. Despite effectiveness, there are several deficiencies in that approach, as it was outlined in \cite{Abdulla1,Abdulla2}:
\begin{itemize}
\item ISP is ill-posed in terms of the phase transition temperature. Moreover, in many applications the latter is not known explicitly, but heavily depends on the process evolution and the environment. For example, in-flight ice accretion on the surface of aircrafts is caused by air droplets which can remain in liquid phase in atmosphere in temperatures much below the freezing temperature, and they impact flying aircraft once they have a surface to freeze on. Similar phenomenon happens during the laser ablation of biomedical tissues. 
\item Iterative gradient descent type methods require solution of the free boundary problem at every iteration, which significantly affects computational cost and accuracy.
\end{itemize} 
To overcome these issues, in \cite{Abdulla1, Abdulla2} a new optimal control framework for the one-phase Stefan problem was developed, where the unknown free boundary is treated as one of the control parameters. The mathematical trick allowed to handle situations with erroneous information on the phase transition temperature, and opened a way to develop numerical methods with reduced computational cost due to the fact that the state vector is a solution of the PDE problem in a fixed region rather than free boundary problem. Frechet differentiability and optimality condition in the new optimal control framework was proved in  \cite{Abdulla3,Abdulla4}, and iterative gradient method for the numerical solution was implemented in \cite{Abdulla5,Abdulla8}.

The approach introduced in \cite{Abdulla1,Abdulla2} is specifically designed for one phase Stefan-type free boundary problems, and is not applicable to multiphase free boundary problems. In a recent paper \cite{Abdulla7}, a new method was introduced for optimal control of multidimensional multiphase Stefan problem based on the weak formulation of the latter as a singular PDE problem with discontinuous coefficient. The idea turned out to be very powerful, allowing to address the above mentioned deficiencies in existing methods, and opened a perspective to develop effective methods for solving the problem under the minimum regularity assumptions on the data. In \cite{Abdulla7} the idea is applied to optimal control of multiphase Stefan problem, where the density of the sources is a control parameter, and the minimization is pursued for the mismatch of the final moment temperature distribution with the desired temperature distribution. Existence of the optimal control and convergence of the finite-difference discretizations of the optimal control problem both with respect to functional and control is proved. In earlier paper \cite{PoggiAbdulla} the method was applied to one dimensional multiphase Stefan problem. In \cite{AbdullaCosgrove} the results of \cite{PoggiAbdulla} are extended to general second order parabolic free boundary problems in space dimension one. The goal of this paper is to extend the method and results of \cite{Abdulla7} to optimal control of singular PDE modeling multiphase Stefan-type free boundary problems for the general second order parabolic operators.

\subsection{Weak Solution of the Singular PDE}\label{OCP}
In this subsection we define the notion of the weak solution of the singular PDE problem \eqref{PDE}-\eqref{vgfuture}. Throughout the paper standard notation of Sobolev spaces will be employed \cite{LSU}.
\begin{definition}\label{typeB} We say that a measurable function $B(x,t,v)$ is \textit{of type }$\n B$ if
    \begin{enumerate}[(a)]
        \item $B(x,t,v)=\beta(v),\qquad v\neq v^j,\quad\forall j=\overline{1,J}$
        \item $B(x,t,v)\in[\beta(v^j)^-,\beta(v^j)^+],\qquad v=v^j$ for some $j$. \\
    \end{enumerate}
    \end{definition}
It should be pointed out that $B(x,t,v)$ can have different values for different $(x,t)$ when $v=v^j$, $j=\overline{1,d}$.

Given $f$, a weak solution $v=v(x,t; f)$ of the problem (\ref{PDE})-(\ref{vgfuture}) is defined as follows:

\begin{definition}\label{weaksoldef} $v\in\overset{\circ}{W}{}_2^{1,1}(D)\cap L_{\infty}(D)$ is called a \emph{weak solution of the singular PDE problem}  (\ref{PDE})-(\ref{vgfuture}) if for some functions $B,B_0$ of type $\n B$, the integral identity
\begin{gather}
\int\limits_D\Big[-B(x,t,v(x,t))\psi_t+\sum\limits_{i=1}^d[a_{i}(x,t)v_{x_i}+b_i(x,t)v]\psi_{x_i}+\sum\limits_{i=1}^d c_i(x,t)v_{x_i}\psi+r(x,t)v\psi-f\psi\Big]\,dxdt  \nonumber \\
- \int\limits_{\Omega}B_0(x,0,\Phi(x))\psi(x,0)\,dx=0\label{weaksol}
\end{gather}
is satisfied for arbitrary $\psi\in\overset{\circ}{W}{}_2^{1,1}(D)$ with { $\psi|_{\Omega\times\{t=T\}}=0$}.
\end{definition}

\subsection{Discrete Optimal Control Problem}\label{DOCP}
We employ a discrete framework introduced in \cite{Abdulla7} to pursue a discretization of the \emph{Problem $\m I$}.
Let $n\in\bb N, \tau:=\frac Tn, h>0$, and slice $\bb R^d\times\bb R$ by the planes 
\[
x_i=k_ih, ~i=1,\ldots,d,\quad t=k_0\tau,\qquad\forall k_{\ell}\in\bb Z,~\ell=0,1,\ldots,d,
\]
so that $\bb R^d\times\bb R$ is split into cells of length $h$ in every $x_i$ direction, and of length $\tau$ in the $t$ direction. We assume the following relation between $h$ and $\tau$;
\begin{equation}\label{htau}
	\frac{h}{\tau} \geq \frac{1+2\sum\limits_{i=1}^{d}\| b_i\|_{L_\infty(D)}}{\bar{b}}
\end{equation}
where $\bar{b}$ is defined in \eqref{bbar}. We will use $\Delta$ as a notation to represent a discretization with steps $(\tau,h)$. We consider a partial ordering on the set of discretizations: we say $\D_1\leq\D_2$ if $\tau_1\leq\tau_2$ and $h_1\leq h_2$. We will denote $t_{\ell}=\tau\ell$ for $\ell=\overline{1,n}$. We will consider two multi-indexes, $\alpha=(k_1,k_2,\ldots,k_d,k_0)$ and $\gamma=(k_1,k_2,\ldots,k_d)$. We can also denote $\alpha=(\gamma,k_0)$, and let $\alpha_i$ be the $i-$th component of $\alpha$, provided $i\in\{1,2,\ldots, d\}$ and $\alpha_0$ is the $d+1-$st component of $\alpha$, and $\gamma_i$ is the $i-$th component of $\gamma$. We can thus represent each elementary cell $C^{\alpha}_{\Delta}$ the rectangular prism $R_{\D}^{\gamma}$ uniquely as
\[
C^{\alpha}_{\Delta}=\Big\{(x,t)\in\bb R^d\times\bb R~\big|~k_ih\leq x_i\leq(k_i+1)h,~i=1,\ldots,d;~~(k_0-1)\tau\leq t\leq k_0\tau\Big\}.
\]
\[
R_{\D}^{\gamma}=\Big\{x\in\bb R^d~\big|~k_ih\leq x_i\leq(k_i+1)h,~i=1,\ldots,d\}.
\]
and a superscript $k$ represents the projection of $R_{\D}^{\gamma}$ onto the hyper-plane $t=k\tau$ of $R^{d+1}$:
\[
R^{\gamma,k}_{\D}=\Big\{(x,t)\in\bb R^d\times\bb R~\big|~k_ih\leq x_i\leq(k_i+1)h,~i=1,\ldots,d;~~t=k\tau\}.
\]
Selecting the collection of prisms and cells contained in $\overline D$ and $\overline\Omega$ respectively, 
\[
\n C^D_{\D}=\Big\{C^{\alpha}_{\D}~|~C_{\D}^{\alpha}\subset\overline D, \ \alpha\in \bb Z^{d+1}\Big\}, \ \ \n R^{\Omega}_{\D}=\Big\{R^{\gamma}_{\D}~|~R_{\D}^{\gamma}\subset\overline{\Omega}, \ \gamma \in \bb Z^d\Big\},
\]
discretized domains are defined as follows:
\[
\Omega_{\D}=\bigcup_{R^{\gamma}_{\D}\in\n R^{\Omega}_{\D}}R^{\gamma}_{\D}~\subset\overline{\Omega},\qquad D_{\D}=\bigcup_{C^{\alpha}_{\D}\in\n C^D_{\D}}C^{\alpha}_{\D}~\subset\overline D.
\]

We define a \emph{natural corner} of a prism $R^{\gamma}_{\D}$ as the vertex with the relatively smallest coordinates in respect to the other vertices. We define the natural corner of the cell $C_{\D}^{(\gamma,k)}$ as the vertex whose spatial coordinates coincide with those of the natural corner of $R^{\gamma}_{\D}$ and whose time coordinate is $k\tau$. Furthermore, each cell and prism will be identified by its natural corner. The lateral boundary of $D_{\D}$ is denoted by $S_{\D}$ and interior sets are defined as $D_{\D}'= (D_{\D}\backslash\partial D_{\D})\cup(\Omega_{\D}\times\{t=T\})$ and $\Omega_{\D}^{'}=\Omega_{\D}\backslash\partial\Omega_{\D}$. Similar notation will be used for cells $C^{\alpha}_{\D}$ and prisms $R_{\D}^{\gamma}$. Next, we introduce the lattice of points

\[
\n L_T=\Big\{(x,t)\in\bb R^d\times\bb R~|~\exists\alpha\in\bb Z^{d+1}\text{ s.t. }x_i=k_ih,~i=1,\ldots,d,~~t=k_0\tau\Big\},
\]
\[
\n L=\Big\{x\in\bb R^d~|~\exists\gamma\in\bb Z^{d}\text{ s.t. }x_i=k_ih,~i=1,\ldots,d\Big\}.
\]

We will denote for notational purposes {$y=(x,t)$, $y_{\alpha}=(k_1h,k_2h,\ldots,k_dh,k_0\tau), x_{\gamma}=(k_1h,k_2h,\ldots,k_dh)$}. Bijections $\alpha\mapsto y_{\alpha}$ and $\gamma\mapsto x_{\gamma}$ will be referred as natural bijections.  

Let X be a set with a natural bijection to a collection of multi-indices, $\gamma$ (or $\alpha$). We will denote $\n A(X)$ as the set of multi-indices corresponding to $X$. For $X\subset\bb R^d$ (or $X\subset\bb R^{d+1}$), we define $\n L(X):=\n L\cap X$. For ease of notation, we will write $\n A(Y)$ in lieu of $\n A(\n L(Y))$ or $\n A(\n L_T(Y))$ and denote $\n A:=\n A(\n R^{\Omega}_{\D})$. The latter means the set of all indices $\gamma$ which are in natural bijection with the natural corners of the prisms in $\Omega_{\D}$. In contrast, $\n A(\Omega_{\D}')$ (or $\n A(\Omega_{\D})$) is the set of indices in natural bijection to the lattice points in the interior of $\Omega_{\D}$ (or in $\Omega_{\D}$). It is obvious that $\n A(\Omega_{\D}')$ is a subset of $\n A$. For ease of notation, in $\Sigma$ and other operations requiring subscripts, expressions like ${\gamma\in\n A(X)}$ will be replaced simply with ${\n A(X)}$. 

Given our data in the appropriate Sobolev or Lebesgue spaces of measurable functions, we define corresponding discrete grid functions via Steklov averages. Whenever it is necessary, we extend all functions to slightly wider region with preservation of the norm. In order to construct discrete version of $ \Phi\in W_2^1(\Omega)$ defined on $\n A(\Omega_{\D})$, we construct an extension of $\Phi$ to $\Omega+B_1(0)$, so that the extension is in $W_2^1(\Omega+B_1(0))$. Such an extension is possible, since $\partial\Omega$ is Lipschitz \cite{LSU}. Let
\begin{equation}\label{phid}
    \phi_{\gamma}=\frac1{h^d}\int\limits_{x_1}^{x_1+h}~\int\limits_{x_2}^{x_2+h}\cdots\int\limits_{x_d}^{x_d+h}\phi(x)\,dx,\qquad\gamma\in\n A(\Omega_{\D}),
\end{equation}
where $\phi$ stands for functions $\Phi, \Gamma$, or initial traces of functions such as $a_i(\cdot,0)$. Let
    \begin{equation}\label{aid}
    g_{\alpha}=\frac1{\tau h^d}\int\limits_{t_{k-1}}^{t_k}~\int\limits_{x_1}^{x_1+h}~\int\limits_{x_2}^{x_2+h}\cdots\int\limits_{x_d}^{x_d+h}g(x,t)\,dx\,dt,\qquad\alpha=(\gamma,k)\in\n A(\n C_{\D}^D), \ k\geq 1
    \end{equation}
where $g$ stands for any of the functions $f,r,a_i,b_i,c_i, i=1,...,d$. 
 
Consider approximation of $\beta(v)$ by the sequence of infinitely differentiable functions 
\begin{equation}\label{bn}
b_{n}(v)=\int_{v-\frac{1}{n}}^{v+\frac{1}{n}}\beta(y)\omega_n(v-y)dy,
\end{equation}
where $\omega_n$ is a standard mollifier defined as
\begin{equation}\label{kernel}
\omega_n(v) =\left\{\begin{matrix}\m C n e^{-\frac{1} {1-n^2v^2}},\quad&|v|\leq\frac{1}{n}\\[2mm] 0,\quad&|v|>\frac{1}{n}\end{matrix}\right.
\end{equation}
and the constant $\m C$ is chosen so that $\int\limits_{\bb R}\omega_1(u)\,du =1$. Since $\beta'(v)$ is piecewise-continuous, we also have
\begin{equation}\label{smoothing_derivative}
b_{n}'(v)=\int_{v-\frac{1}{n}}^{v+\frac{1}{n}}\beta'(y)\omega_n(v-y)dy.
\end{equation}
This implies $b_n$ is also strict monotonically increasing and by \eqref{bbar} we have
\begin{equation}\label{b_n bound}
b_n'(v)\geq \bar{b}>0
\end{equation}  
For a given discretization $\D$, we define a finite-dimensional discrete control vector 
\[ [f]_{\D}:=\{f_{\alpha}: \ f_\alpha\in\mathbb{R}, \alpha\in\n A(\n C_{\D}^D)\}\]
and corresponding discrete norms:
\begin{equation}\nonumber
\Vert[f]_{\D}\Vert_{\ell_{\infty}}:=\max\limits_{\n A(\n C_{\D}^D)}|f_{\alpha}|, \ \   \Vert[f]_{\D}\Vert_{\ell_2}:=\Big (\sum\limits_{\n A(\n C_{\D}^D)}\tau h^d f_{\alpha}^2\Big )^{\frac{1}{2}}.
\end{equation}
For any collection $\{v_{\alpha}\}$, with $\alpha=(\gamma, k_0)$, we utilize the notation
\[
\gamma\pm e_i:=(k_1,\ldots,k_i\pm1,\ldots,k_d),\quad \alpha\pm e_i:=(k_1,\ldots,k_i\pm 1,\ldots,k_d,k_0)
\]
for suitable $i$. We employ the standard notation for the first order backward and forward space and time differences:
\[
v_{\alpha\bar t}=\frac{v_{(\gamma,k_0)}-v_{(\gamma,k_0-1)}}{\tau}, \ v_{\alpha t}=\frac{v_{(\gamma,k_0+1)}-v_{(\gamma,k_0)}}{\tau}, \ v_{\alpha x_i}=\frac{v_{\alpha+e_i}-v_{\alpha}}{h}, \ v_{\alpha \bar x_i}=\frac{v_{\alpha}-v_{\alpha-e_i}}{h}.
\]
For a given $R>0$, let 
\[
\n F_{\D}^R:=\Big\{[f]_{\D}~\big|~\Vert[f]_{\D}\Vert_{\ell_{\infty}}\leq R\Big\}
\]
be a \emph{discrete control set}. We define
\[
\n P_{\D}:\bigcup_R\n F_{\D}^R\longrightarrow\bigcup_R\n F^R,\qquad\n P_{\D}([f]_{\D})=f^{\D}
\]
to be the interpolating map from the discrete control set to the continuous control set, where
\[
f^{\D}\Big|_{C_{\D}^{\alpha}}=f_{\alpha},~~\alpha\in\n A(\n C_{\D}^D),\qquad f^{\D}\equiv0~~\text{elsewhere on }D.
\]
Similarly, we define
\[
\n Q_{\D}:\bigcup_R\n F^R\longrightarrow\bigcup_R\n F_{\D}^R,\qquad\n Q_{\D}(f)=[f]_{\D}
\]
to be the discretizing map from the continuous control set to the discrete control set, where $f_{\alpha}$ is given by (\ref{aid}) for each $\alpha\in\n A(\n C_{\D}^D)$.

Next, we are going to define a solution of the discrete singular PDE problem.\\
\begin{definition}\label{dsvdef} Given $[f]_{\D}$, the vector function $[v([f]_{\D})]_{\D}=(v(0),v(1),\ldots,v(n))$, where $v(k)$ is a collection of real numbers $\{v_{\gamma}(k)\}, \gamma\in\n A(\Omega_{\D}),~k=0,1,2\ldots,n$, is called a \emph{discrete state vector} if
    \begin{enumerate}[{\bf(i)}]
    \item $v_{\gamma}(0)=\Phi_{\gamma},~\gamma\in\n A(\Omega_{\D}')$,
    \item For each fixed $k=1,\ldots,n$, the collection $v(k)$ satisfies
    \begin{gather}
    \sum\limits_{{\n A}}h^d\Bigg[\Big(b_n(v_{\gamma}(k))\Big)_{\bar t}\eta_{\gamma}~+\sum\limits_{i=1}^{d}[(a_{i})_\alpha v_{\gamma x_i}(k)+(b_i)_\alpha v_\gamma(k)]\eta_{\gamma x_i}~+\sum\limits_{i=1}^{d} (c_i)_\alpha v_{\gamma x_i}(k) \eta_\gamma\nonumber \\
     ~+r_\alpha v_\gamma(k)\eta_\gamma ~-f^{\D}_{(\gamma,k)}\eta_{\gamma}\Bigg]=0\label{dsveq}
    \end{gather}
    for arbitrary $\{\eta_{\gamma}\},~\gamma\in\n A(\Omega_{\D})$ such that $\eta_{\gamma}=0$ for $\gamma\in\n A(\partial\Omega_{\D})$.
    \item For each $k=0,1,\ldots,n$, we have $v_{\gamma}(k)=0$ for $\gamma\in\n A(\partial\Omega_{\D})$.
    \end{enumerate}
    \end{definition}
    
    It should be mentioned that the collection $\{f^{\D}_{\alpha}\}$ in (\ref{dsveq}) coincides with $\n Q_{\D}\Big(\n P_{\D}([f]_{\alpha})\Big)$. It will be proved in Lemma~\ref{exuniqdsv}, Section~\ref{prelres} that for any $[f]_{\D}\in\n F_{\D}^R$ there exists a unique discrete state vector. We can thus define the discrete cost functional $\n I_{\D}:\cup_R\n F_{\D}^R\ra[0,+\infty)$ by
    \begin{equation}\label{discretecostfunctional}
    \n I_{\D}([f]_{\D})=\sum\limits_{\n A}h^d|v_{\gamma}(n)-\Gamma_{\gamma}|^2
    \end{equation}
    where $v_{\gamma}(n)$ represents the $n$th component of the discrete state vector $[v([f]_{\D})]_{\D}$. We will refer to the discrete optimal control problem on the minimization of the functional \eqref{discretecostfunctional} on a discrete control set $\n F_{\D}^R$ as a \emph{Problem $\m I_{\D}$}. 
    
Next, we introduce various interpolations of the discrete state vector $[v]_\D$. First, define a piecewise constant interpolation $\tilde V_{\D}:D\ra\bb R$, which assigns the value of $[v]_\D$ on the natural corner to the interior and top face of each cell, i.e.
\begin{equation}\label{pwconstant}
    \tilde V_{\D}\Big|_{C_{\D}^{\alpha'}\cup R_{\D}^{\gamma,k}}=v_{\gamma}(k),\qquad\forall\alpha=(\gamma,k)\in\n A(\n C_{\D}^D),
\end{equation}
with $\tilde V_{\D}=0$ everywhere else in $D$. We also define a piecewise constant interpolation of the discrete $x_i$-derivative denoted as $\tilde V^i_{\D}:D\ra\bb R$, which assigns the value of the forward $x_i$-difference of $[v]_\D$ at the natural corner to interior and top face of each cell, i.e.
\begin{equation}\label{derivativestep}
    \tilde V_{\D}^i\Big|_{C_{\D}^{\alpha'}\cup R_{\D}^{\gamma,k}}=v_{\gamma x_i}(k),\qquad\forall\alpha=(\gamma,k)\in\n A(\n C_{\D}^D)
\end{equation}
with $\tilde V^i_{\D}=0$ everywhere else in $D$. 

For fixed $k=\overline{0,n}$, we define a multilinear interpolation $V_{\D}^k:\Omega\ra\bb R$ as a function that takes the value $v_\gamma(k)$ at corresponding lattice points of $\Omega_\D$, is linear with respect to every variable, when all other variables are fixed, and vanishes in $\Omega\backslash\Omega_\D$. Note that $V_\D^k\in C(\bar{\Omega})\cap W_2^1(\Omega)$. Next, we define $V_{\D}:D\ra\bb R$ as a piecewise constant interpolation of $V_{\D}^k$ onto $[0, T]$:
\begin{equation}\label{pwlinearpwconstant}
    V_{\D}(x,t)=V_{\D}^k(x),\qquad t\in(t_{k-1},t_k],~k=1,2,\ldots,n
\end{equation}
with $V_{\D}(x,0)=V_{\D}^0(x)$. We have $V_\D\in W_2^{1,0}(D)$. Finally, we define multilinear interpolation  $V_\D'\in W_2^{1,1}(D)\cap C(\bar{D})$:
\begin{equation}\label{pwlinear}
    V_{\D}'(x,t)=V_{\D}^{k-1}(x)+\big(V_{\D}^k(x)\big)_{\bar t}(t-t_{k-1}),\qquad t\in[t_{k-1},t_k],~k=1,2,\ldots,n.
\end{equation}

\subsection{Main Results}\label{mainresults}
Throughout the paper we assume the following assumptions are satisfied:
\begin{gather}
    \Phi\in W_2^1(\Omega) \cap L_\infty(\Omega), \quad \Gamma\in L_2(\Omega)\label{phibound}\\
    a_i, b_i, c_i, \in W^{1,0}_{\infty}(D), \quad \frac{\partial a_i}{\partial t} \in L_{\infty,1}(D), \ \  i = \overline{1,d}, \quad r \in L_\infty(D), \label{databounds}\\
 \{x\in\Omega\, |\, \Phi(x)=v^j\}, j=\overline{1,m} \ \text{has d-dimensional Lebesgue measure 0},\label{phiset}
 \end{gather}
$a_i, i=1,...,d$ satisfy \eqref{apositive}; \ \ $\beta$ is a maximal monotone graph satisfying \eqref{beta},\eqref{bbar}.

The following are the main results of this paper:

\begin{theorem}\label{optsol} There exists an optimal control in problem $\m I$, i.e.
	\[
	\n F_*:=\Big\{f\in\n F^R~\Big|~\n J(f)=\n J_*:=\inf\limits_{f\in\n F^R}\n J(f)\Big\}\neq \emptyset
	\]
\end{theorem}

\begin{theorem}\label{funcapprox}The sequence of discrete optimal control problems $\m I_n$ approximates the optimal control problem $\m I$ with respect to the functional and control, i.e.    \begin{equation}\label{approx1}
    \lim\limits_{\D\ra0}\n I_{\D_*}=\n J_*
    \end{equation}
    where
    \[
    \n I_{\D_*}=\inf\limits_{\n F_{\D}^R}\n I([f]_{\D});
    \]
    if $[f]_{\D,\ep}\in\n F_{\D}^R$ is chosen such that 
        \begin{equation}\label{approxcond}
    \n I_{\D_*}\leq\n I_{\D}([f]_{\D,\ep})\leq\n I_{\D_*}+\ep_{\D}, \ \  \ep_{\D}\downarrow 0,
    \end{equation}
    then we have
    \begin{equation}\label{approx2}
    \lim\limits_{\D\ra0}\n J(\n P_{\D}([f]_{\D,\ep}))=\n J_*,
    \end{equation}
the sequence $\{\n P_{\D}([f]_{\D,\ep})\}$ is weakly precompact in $L_2(D)$, and all of its weak limit points lie in $\n F_*$. Moreover, if $f_*$ is such a weak limit point, then there is a subsequence $\D'$ such that the multilinear interpolations $V_{\D'}'$ of the discrete state vectors $[v([f]_{\D',\ep})]_{\D'}$ converge to the weak solution $v=v(x,t;f_*)\in W_2^{1,1}(D)\cap L_\infty(D)$ of the singular PDE problem \eqref{PDE}-\eqref{vgfuture}, weakly in $W_2^{1,1}(D)$, strongly in $L_2(D)$, and almost everywhere on $D$.     \end{theorem}

\section{Preliminary Results}\label{prelres}
We define
    \begin{equation}\label{zeta1}
    \zeta_{\D}^{\gamma,k}:=\int\limits_0^1b_n'\Big(\theta v_{\gamma}(k)+(1-\theta)v_{\gamma}(k-1)\Big)\,d\theta,
    \end{equation}
    for each $(\gamma,k)\in\n A(D_{\D}), k\neq0$. Note that for every $(\gamma,k)$, we have
    \begin{equation}\label{mvt}
    \Big(b_n(v_{\gamma}(k))\Big)_{\bar t}=\zeta_{\D}^{\gamma,k}v_{\gamma\bar t}(k).
    \end{equation}
    and
\begin{equation}\label{zetabig}
\zeta_{\D}^{\gamma,k}\geq\inf\limits_{x\in\bb R}b_n'(x)\geq\bar b,
\end{equation}
independently of $\gamma, k,\D, v_\gamma(k)$.

\begin{lemma}\label{equiv} For a fixed discretization $\D$ and given discrete control $[f]_{\D}$, a vector function $[v([f]_{\D})]_{\D}$ is a discrete state vector in the sense of Definition~\ref{dsvdef} if and only if it satisfies conditions (i), (ii)', and (iii), where
	
    \textbf{(ii)'}$~$ $\forall k=\overline{1,n}$ and $\gamma\in\n A(\Omega_{\D}')$, we have
    \begin{gather}
    \Big(b_n(v_{\gamma}(k))\Big)_{\bar t}-\sum\limits_{i=1}^{d}\Big([(a_{i})_\alpha v_{\gamma x_i}(k)+(b_i)_\alpha v_\gamma(k)]\Big)_{\bar{x}_i} \nonumber \\
    ~+\sum\limits_{i=1}^{d} (c_i)_\alpha v_{\gamma x_i}(k) ~+r_\alpha v_\gamma(k)=f^{\D}_{(\gamma,k)}. \label{dsvequiv}
    \end{gather}
\end{lemma}

\emph{Proof.} Assume $[v([f]_{\D})]_{\D}$ satisfy (i),(ii)' and (iii), and $k\in\{1,...,n\}$ is fixed. Take $\{\eta_{\gamma}\in\mathbb{R}: \ \gamma\in\n A(\Omega_{\D})\}$ with $\eta_{\gamma}=0$ for $\gamma\in\n A(\partial\Omega_{\D})$. Multiplying (\ref{dsvequiv}) by $h^d\eta_{\gamma}$, and performing summation with respect to $\gamma\in\n A(\Omega_{\D}')$ we have
\begin{align}\label{almostdsveq}
\sum\limits_{\n A(\Omega'_{\D})}h^d&\Bigg[\Big(b_n(v_{\gamma}(k))\Big)_{\bar t}\eta_{\gamma}~-\sum\limits_{i=1}^{d}\Big([(a_{i})_\alpha v_{\gamma x_i}(k)+(b_i)_\alpha v_\gamma(k)]\Big)_{\bar{x}_i}\eta_{\gamma}\nonumber\\
&~+\sum\limits_{i=1}^{d} (c_i)_\alpha v_{\gamma x_i}(k)\eta_{\gamma} ~+r_\alpha v_\gamma(k)\eta_{\gamma}-f^{\D}_{(\gamma,k)}\eta_{\gamma}\Bigg]=0.
\end{align}
For any fixed $i$, we can rewrite
\begin{gather*}
-\sum\limits_{\n A(\Omega'_{\D})}\Big([(a_{i})_\alpha v_{\gamma x_i}(k)+(b_i)_\alpha v_\gamma(k)]\Big)_{\bar{x}_i}\eta_{\gamma}\\
=-\sum\limits_{\n A(\Omega'_{\D})}\frac{[(a_{i})_\alpha v_{\gamma x_i}(k)+(b_i)_\alpha v_\gamma(k)]}h\eta_{\gamma}+\sum\limits_{\n A(\Omega'_{\D})}\frac{[(a_{i})_{\alpha-e_i} v_{(\gamma-e_i) x_i}(k)+(b_i)_{\alpha-e_i} v_{\gamma-e_i}(k)]}h\eta_{\gamma},\\
=-\sum\limits_{\n A(\Omega'_{\D})}\frac{[(a_{i})_\alpha v_{\gamma x_i}(k)+(b_i)_\alpha v_\gamma(k)]}h\eta_{\gamma}
+\sum\limits_{\gamma\text{ s.t. } \gamma+e_i\in\n A(\Omega'_{\D})}\frac{[(a_{i})_{\alpha} v_{\gamma x_i}(k)+(b_i)_{\alpha} v_\gamma(k)]}h\eta_{\gamma+e_i}
\\=\sum\limits_{\gamma\text{ s.t. } \gamma\in\n A(\Omega_{\D}')\text{ and }\gamma+e_i\in\n A(\Omega'_{\D})}[(a_{i})_\alpha v_{\gamma x_i}(k)+(b_i)_\alpha v_\gamma(k)]\eta_{\gamma x_i}\\
-\sum\limits_{\gamma\text{ s.t. } \gamma\in\n A(\Omega_{\D}')\text{ and }\gamma+e_i\in\n A(\partial\Omega_{\D})}\frac{[(a_{i})_\alpha v_{\gamma x_i}(k)+(b_i)_\alpha v_\gamma(k)]}h\eta_{\gamma}~\\
+\sum\limits_{\gamma\text{ s.t. }\gamma\in\n A(\partial\Omega_{\D})\text{ and }\gamma+e_i\in\n A(\Omega'_{\D})}\frac{[(a_{i})_{\alpha} v_{\gamma x_i}(k)+(b_i)_{\alpha} v_\gamma(k)]}h\eta_{\alpha+e_i}\\
=\sum\limits_{\gamma\text{ s.t. }\gamma\in\n A(\Omega_{\D}')\text{ and } \gamma+e_i\in\n A(\Omega'_{\D})}[(a_{i})_\alpha v_{\gamma x_i}(k)+(b_i)_\alpha v_\gamma(k)]\eta_{\gamma x_i}
\end{gather*}
\begin{gather*}
+\sum\limits_{\gamma\text{ s.t. }\gamma\in\n A(\Omega_{\D}')\text{ and }\gamma+e_i\in\n A(\partial\Omega_{\D})}\frac{[(a_{i})_\alpha v_{\gamma x_i}(k)+(b_i)_\alpha v_\gamma(k)]}h\big(-\eta_{\gamma}+\eta_{\gamma+e_i}\big)~\\
+\sum\limits_{\gamma\text{ s.t. }\gamma\in\n A(\partial\Omega_{\D})\text{ and }\gamma+e_i\in\n A(\Omega'_{\D})}\frac{[(a_{i})_{\alpha} v_{\gamma x_i}(k)+(b_i)_{\alpha} v_\gamma(k)]}h\big(\eta_{\gamma+e_i}-\eta_\gamma\big)\\
 =\sum\limits_{\n A} [(a_{i})_\alpha v_{\gamma x_i}(k)+(b_i)_\alpha v_\gamma(k)]\eta_{\gamma x_i}.
\end{gather*}
Taking into account this transformation in (\ref{almostdsveq}), (ii) easily follows. Now conversely, suppose (i), (ii) and (iii) are satisfied. We fix any $k\in\overline{1,n}$, and any $\gamma'\in\n A(\Omega'_{\D})$, and choose the collection $\{\eta_{\gamma}\}$ as follows: $\eta_{\gamma'}=1$, and $\eta_{\gamma}=0,\, \forall \gamma \neq \gamma'$. Then from (\ref{dsveq}) we get

\begin{gather*}
    \Big(b_n(v_{\gamma'}(k))\Big)_{\bar t}+\sum\limits_{i=1}^d\left(-\frac{[(a_{i})_{\alpha'} v_{\gamma' x_i}(k)+(b_i)_{\alpha'} v_{\gamma'}(k)]}h\right )\\
    \sum\limits_{i=1}^d\frac{[(a_{i})_{\alpha'-e_i} v_{(\gamma'-e_i) x_i}(k)+(b_i)_{\alpha'-e_i} v_{\gamma'-e_i}(k)]}h
    +\sum\limits_{i=1}^{d} (c_i)_{\alpha'} v_{\gamma' x_i}(k) ~+r_{\alpha'} v_{\gamma'}(k)-f^{\D}_{(\gamma',k)}=0
\end{gather*}
which yields (\ref{dsvequiv}) for $\gamma'$. Since $\gamma' \in \n A(\Omega'_{\D})$ is arbitrary, statement (ii)' follows.\hfill{$\square$}

\begin{lemma}\label{exuniqdsv} For any discretization $\D$ with sufficiently small $h$ and $\tau$ satisfying \eqref{htau}, and for any  $[f]_{\D}\in\n F^R_{\D}$, there exists a unique discrete state vector $[v([f]_{\D})]_{\D}$.
\end{lemma}
\emph{Proof.} To prove uniqueness, assume that $[v([f_{\D}])]_{\D}, [\tilde v([f_{\D}])]_{\D}$ are two discrete state vectors. We use induction on $k$. We have $v(0)=\tilde{v}(0)$, due to conditions (i) and (iii). Fix any $k,~1\leq k\leq n$ and assume $v(k-1)=\tilde v(k-1)$. By selecting $\eta=v(k)-\tilde v(k)$, and by subtracting \eqref{dsveq} for both $v(k)$ and $\tilde v(k)$, we derive
\begin{gather*}
    \sum\limits_{\n A}h^d\Bigg[\Big((b_n(v_{\gamma}(k)))_{\bar t}-(b_n(\tilde v_{\gamma}(k)))_{\bar t}\Big)(v_{\gamma}(k)-\tilde v_{\gamma}(k))+\sum\limits_{i=1}^d\Big[ (a_{i})_\alpha(v_{\gamma x_i}(k)-\tilde{v}_{\gamma x_i}(k))^2\Big]\\
    +\sum\limits_{i=1}^d\Big[ \Big((b_i)_\alpha+(c_i)_\alpha)(v_\gamma(k)-\tilde{v}(k))(v_{\gamma x_i}(k)-\tilde{v}_{\gamma x_i}(k))\Big]+r_\alpha(v_\gamma(k)-\tilde{v}_\gamma(k))^2 \Bigg]=0,
\end{gather*}
We can rewrite the following:
\begin{align*}
(b_n(v_{\gamma}(k)))_{\bar t}-(b_n(\tilde v_{\gamma}(k)))_{\bar t}&=
\frac{b_n(v_{\gamma}(k))-b_n(\tilde v_{\gamma}(k))}{\tau}
\end{align*}
Using \eqref{apositive}, \eqref{bbar}, the integral mean value theorem as in \eqref{mvt},\eqref{zetabig}, and Cauchy inequality with $\epsilon >0$, we get
\begin{gather*}
	\sum\limits_{\n A}h^d\Bigg[\frac{\bar{b}}{\tau}(v_{\gamma}(k)-\tilde v_{\gamma}(k))^2+a_0 \sum\limits_{i=1}^{d}(v_{\gamma x_i}(k)-\tilde{v}_{\gamma x_i}(k))^2+r_\alpha(v_\gamma(k)-\tilde{v}_\gamma(k))^2\Big]\leq\\
	 \frac{a_0}{2}\sum\limits_{\n A} h^d\sum\limits_{i=1}^d \Big[(v_{\gamma x_i}(k)-\tilde{v}_{\gamma x_i}(k))^2 \Big] + \frac{1}{a_0}\sum\limits_{\n A}h^d \Bigg[\sum\limits_{i=1}^d(\|b_i\|_{L_\infty}^2+\|c_i\|_{L_\infty}^2)(v_\gamma(k)-\tilde{v}_\gamma(k))^2  \Bigg],
\end{gather*}
and therefore,
\begin{gather*}
	\sum\limits_{\n A}h^d\Bigg[\Bigg(\frac{\bar{b}}{\tau} - \frac{1}{a_0}\sum\limits_{i=1}^d(\|b_i\|_{L_\infty}^2+\|c_i\|_{L_\infty}^2) - \|r\|_{L_\infty}\Bigg)(v_{\gamma}(k)-\tilde v_{\gamma}(k))^2\nonumber\\ 
	+ \frac{a_0}{2}\sum\limits_{i=1}^{d}(v_{\gamma x_i}(k)-\tilde{v}_{\gamma x_i}(k))^2\Bigg] \leq 0
\end{gather*}

By taking $\tau$ sufficiently small, all the terms on the left hand side become non-negative, and therefore, each term is equal to 0. This implies that $v_{\gamma}(k)=\tilde v_{\gamma}(k)$ for $\gamma\in\n A(\Omega'_{\D})$. By (iii) and induction argument, we get $v=\tilde v$, and uniqueness follows.\\

Now we prove existence, again through induction on $k$. Let discretization $\D$ and $[f]_{\D}$ are fixed. For $k=0$, $v(0)$ is given by (i) and (iii) of Definition \ref{dsvdef}. Assuming that $v(0),v(1),\ldots,v(k-1)$ exist, we prove the existence of $v(k)$ by the method of successive approximations. By (iii), $v(k)$ taken to be $0$ for any lattice point on the boundary of $\Omega_{\D}$. For the interior lattice points, we rewrite (\ref{dsvequiv}) as

\begin{gather}
    \frac{h^2}{\tau}\big[b_n(v_{\gamma}(k))-b_n(v_{\gamma}(k-1))\big]+\Bigg[\sum\limits_{i=1}^{d}\Big((a_i)_\alpha+(a_i)_{\alpha-e_i}-h(b_i)_\alpha-h(c_i)_\alpha\Big)+h^2r_\alpha \Bigg]v_\gamma(k) \nonumber \\
    -\sum\limits_{i=1}^{d}\Bigg[ \Big((a_i)_\alpha-h(c_i)_\alpha\Big)v_{\gamma+e_i}(k)+\Big((a_i)_{\alpha-e_i}-h(b_i)_{\alpha-e_i}\Big)v_{\gamma-e_i}(k)\Bigg]=h^2f^{\D}_{(\gamma,k)}. \label{vsystem}
\end{gather}
We set $v^0=v(k-1)$, and having calculated $v^N$, $v^{N+1}$ is found as a solution of the system
\begin{gather}\label{vNsystem}
\frac{h^2}{\tau}b_n(v_{\gamma}^{N+1})+\Bigg[\sum\limits_{i=1}^{d}\Big((a_i)_\alpha+(a_i)_{\alpha-e_i}-h(b_i)_\alpha-h(c_i)_\alpha\Big)+h^2r_\alpha \Bigg]v_{\gamma}^{N+1}
=\frac{h^2}{\tau}b_n(v_{\gamma}(k-1))\nonumber\\
+\sum\limits_{i=1}^{d}\Bigg[ \Big((a_i)_\alpha-h(c_i)_\alpha\Big)v_{\gamma+e_i}^N+\Big((a_i)_{\alpha-e_i}-h(b_i)_{\alpha-e_i}\Big)v_{\gamma-e_i}^N\Bigg]+h^2f^{\D}_{(\gamma,k)}. \label{vNsystem}
\end{gather}

Since the left hand side of \eqref{vNsystem} is monotonically increasing with respect to $v^{N+1}$ for sufficiently small $h$, and has a range $\mathbb{R}$, there exists a unique solution $v^{N+1}$. This implies the sequence $\{v^N\}$ is well-defined. Subtracting (\ref{vNsystem}) for $N$ and $N-1$ we have
\begin{gather}
    \frac{h^2}{\tau}\Big(b_n(v_{\gamma}^{N+1})-b_n(v_{\gamma}^N)\Big)+\Bigg[\sum\limits_{i=1}^{d}\Big((a_i)_\alpha+(a_i)_{\alpha-e_i}-h(b_i)_\alpha-h(c_i)_\alpha\Big)+h^2r_\alpha \Bigg]\Big(v_{\gamma}^{N+1}-v_{\gamma}^N\Big) \nonumber \\
    = \sum\limits_{i=1}^{d}\Bigg[ \Big((a_i)_\alpha-h(c_i)_\alpha\Big)(v_{\gamma+e_i}^N-v_{\gamma+e_i}^{N-1})+\Big((a_i)_{\alpha-e_i}-h(b_i)_{\alpha-e_i}\Big)(v_{\gamma-e_i}^{N}-v_{\gamma-e_i}^{N-1})\Bigg]\label{system2}
\end{gather}
Similar to \eqref{mvt},\eqref{zeta1}, we have
\begin{equation}\label{zeta}
b_n(v_{\gamma}^{N+1})-b_n(v_{\gamma}^N)=\zeta_{\D,N}^{\gamma,k}\Big(v_{\gamma}^{N+1}-v_{\gamma}^N\Big), \ \zeta_{\D,N}^{\gamma,k}:=\int\limits_0^1b_n'\Big(\theta v_{\gamma}^{N+1}+(1-\theta)v_{\gamma}^N\Big)\,d\theta,
\end{equation}
where $\zeta_{\D,N}^{\gamma,k}$ satisfies (\ref{zetabig}) uniformly with respect to $\D, \gamma, k, N$. From (\ref{system2}), \eqref{zeta} it follows that
\begin{gather}
    v_{\gamma}^{N+1}-v_{\gamma}^N =\frac1{\frac{h^2}{\tau}\zeta_{\D,N}^{\gamma,k}+\sum\limits_{i=1}^{d}\Big((a_i)_\alpha+(a_i)_{\alpha-e_i}-h(b_i)_\alpha-h(c_i)_\alpha\Big)+h^2r_\alpha} \times \nonumber \\
    \Bigg(\sum\limits_{i=1}^{d}\Bigg[ \Big((a_i)_\alpha-h(c_i)_\alpha\Big)(v_{\gamma+e_i}^N-v_{\gamma+e_i}^{N-1})+\Big((a_i)_{\alpha-e_i}-h(b_i)_{\alpha-e_i}\Big)(v_{\gamma-e_i}^{N}-v_{\gamma-e_i}^{N-1})\Bigg]\Bigg)\label{system3}
\end{gather}

Due to (\ref{zetabig}),\eqref{apositive},\eqref{databounds}, for sufficiently small $h$ we have 
\begin{gather*}
	0<\frac1{\frac{h^2}{\tau}\zeta_{\D,N}^{\gamma}+\sum\limits_{i=1}^{d}\Big((a_i)_\alpha+(a_i)_{\alpha-e_i}-h(b_i)_\alpha-h(c_i)_\alpha\Big)+h^2r_\alpha}\\
	 \leq \frac1{\frac{h^2}{\tau}\bar{b}+\sum\limits_{i=1}^{d}\Big((a_i)_\alpha+(a_i)_{\alpha-e_i}-h(b_i)_\alpha-h(c_i)_\alpha\Big)+h^2r_\alpha}
\end{gather*}
Let
\[
A_N:=\max\limits_{\gamma}|v_{\gamma}^{N+1}-v_{\gamma}^N|.
\]
From \eqref{system3} we deduce that for sufficiently small $h$ and for every $\gamma$
\begin{equation}\label{A_Ninequality}
|v_{\gamma}^{N+1}-v_{\gamma}^N|\leq\delta A_{N-1},
\end{equation}
where
\begin{gather*}
\delta:=
\Biggr(1+\frac{\frac{h^2}{\tau}\bar{b}	+h\sum\limits_{i=1}^{d}(b_i)_{\alpha-e_i}-h\sum\limits_{i=1}^{d}(b_i)_{\alpha}+h^2r_\alpha}{\sum\limits_{i=1}^{d}\Big( (a_i)_\alpha+(a_i)_{\alpha-e_i}-h(b_i)_{\alpha-e_i}-h(c_i)_\alpha\Big)}\Biggr)^{-1}
\end{gather*}
From \eqref{apositive},\eqref{databounds},\eqref{htau} it follows that for sufficiently small $h$
\begin{gather*}
\sum\limits_{i=1}^{d}\Big( (a_i)_\alpha+(a_i)_{\alpha-e_i}-h(b_i)_{\alpha-e_i}-h(c_i)_\alpha\Big)\geq 2da_0-h\sum\limits_{i=1}^d\Big (\|b_i\|_{L_\infty(D)}+(\|c_i\|_{L_\infty(D)}\Big ) >0,\nonumber\\
	\frac{h}{\tau}\bar{b}+\sum\limits_{i=1}^{d}(b_i)_{\alpha-e_i}-\sum\limits_{i=1}^{d}(b_i)_{\alpha}+hr_\alpha \geq \frac{h}{\tau}\bar{b}-2\sum\limits_{i=1}^{d}\| b_i\|_{L_\infty(D)}-h\|r\|_{L_\infty(D)}\geq 1-h \|r\|_{L_\infty(D)}> 0.
\end{gather*}
Hence, $\delta \in (0,1)$, and by taking maximum with respect to $\gamma$ from \eqref{A_Ninequality} we derive inductive chain of inequalities
\begin{equation}\label{chain}
    A_N\leq\delta A_{N-1}\leq\delta^2 A_{N-2}\leq\cdots\leq\delta^N A_0.
\end{equation}
Following the proof of the Lemma~7, \cite{Abdulla7}, from \eqref{chain} it follows that there exists a limit
\begin{equation}\label{successive}
v_{\gamma}(k)=\lim\limits_{N\ra\infty}v_{\gamma}^N,\quad\gamma\in\n A(\Omega_{\D}').
\end{equation}
and $v(k)$, given by (\ref{successive}) satisfies (\ref{dsveq}). Thus, the existence of the discrete state vector is proved.\hfill{$\square$}

The uniqueness of the weak solution of the multiphase Stefan problem, or singular PDE problem \eqref{PDE}-\eqref{vgfuture} with $\mathcal{L}=\Delta$, in the sense of Definition~\ref{weaksoldef}, is proved in \cite{LSU}. Next proposition formulates uniqueness of the weak solution of the singular PDE problem \eqref{PDE}-\eqref{vgfuture}. 
\begin{proposition}\label{some} There exists at most one solution $v\in\overset{\circ}{W}{}_2^{1,1}(D)\cap L_{\infty}(D)$ of the singular PDE problem \eqref{PDE}-\eqref{vgfuture}. \end{proposition}
\emph{Proof.} We prove uniqueness in a broader class of solutions $v\in L_{\infty}(D)$, which satisfy the following integral identity instead of \eqref{weaksol}:
\begin{equation}\label{weakersol}
    \int\limits_D\Big[B(x,t,v)\psi_t+v \mathcal{L}^*\psi+f\psi\Big]\,dxdt +\int\limits_\Omega B_0(x,0,\Phi(x))\psi(x,0)\,dx =0, 
    \end{equation}
$\forall \psi\in W_2^{2,1}(D)$ such that $\psi(x,T)\vert_{x \in \Omega}=0,  \psi \vert_{\partial \Omega \times (0,T]}=0$, where
\[ \mathcal{L}^*\psi=\sum\limits_{i=1}^d(a_i\psi_{x_i})_{x_i}-\sum\limits_{i=1}^d b_i \psi_{x_i}+\sum\limits_{i=1}^d(c_i\psi)_{x_i}-r\psi.\]
Subtracting any solutions $v, \tilde{v}\in L_{\infty}(D)$ of \eqref{weakersol}, and by taking into account \eqref{phiset}, we get
\begin{equation}\label{weakersol1}
\int\limits_D \big(B(x,t,v)-\tilde B(x,t,\tilde v)\big)\left(\psi_t+z(x,t)\mathcal{L}^*\psi \right)\,dx\,dt = 0,
\end{equation}
where 
\begin{equation*} 
z(x,t)= 
\left\{
    \begin{array}{l}
    \frac{v-\tilde v}{B(x,t,v)-\tilde B(x,t,\tilde v)}, \ \text{if} \ v(x,t)\neq \tilde{v}(x,t),\\
    0, \ \text{if} \ v(x,t)= \tilde{v}(x,t).
    \end{array}\right.
\end{equation*}
Since $B,\tilde{B} \in \n B$, we have 
\begin{equation}\label{zbound}
0\leq z \leq \frac{1}{\bar{b}}, \ \text{a.e.} \ (x,t)\in D. 
\end{equation}
Fix $\ep>0$, and take $\psi(x,t)$ to be the solution of the Dirichlet problem for the backward parabolic PDE:
\begin{gather}
\psi_t+z^\ep(x,t)\mathcal{L}^*\psi=F(x,t), \ \text{in} \ \Omega\times [0,T)\label{conjugate},
\\
\psi(x,T)\vert_{x \in \Omega}=0, \ \  \psi \vert_{\partial \Omega \times (0,T]}=0, \label{psiboundary}
\end{gather}
where $z^\ep(x,t)=z(x,t)+\ep$, $F$ is an arbitrary compactly-supported, smooth function in $D$. From \cite{LSU} it follows there exists a unique solution $\psi^{\ep} \in W_2^{2,1}(D)$.  By using \eqref{weakersol1},\eqref{conjugate}, we can write
\begin{equation}
   \int\limits_D\hat{B}(x,t)\left(F-\ep \mathcal{L}^*\psi  \right)\,dx\,dt = 0 \label{almost}.
\end{equation}
where $\hat{B}(x,t)=B(x,t,v(x,t))-\tilde B(x,t,\tilde v(x,t))$. Our goal is to eliminate the $\ep$-term by passing to limit as $\ep \downarrow 0$, and use the arbitrariness of $F$ to derive that $\hat{B}=0$ a.e. on $D$. 
To do that we need to attain energy estimate for the solution of \eqref{conjugate},\eqref{psiboundary}. For simplicity, we will derive the required energy estimate for the parabolic PDE by assuming that time variable $t$ is replaced with $T-t$ in \eqref{conjugate},\eqref{psiboundary}. Let us multiply the parabolic version of (\ref{conjugate}) by $\sum\limits_{i=1}^d(a_i\psi_{x_i})_{x_i}$, integrate it over $D_t :=\Omega\times(0,t)$, to get
\begin{equation}
    -\int\limits_{D_t} (\psi_{\tau}-z^\ep\mathcal{L}^*\psi)\sum\limits_{i=1}^{d}(a_i \psi_{x_i})_{x_i}\,dx\,d\tau = \int\limits_{D_t} \sum\limits_{i=1}^da_i\psi_{x_i}F_{x_i}\,dx\,d\tau, \label{psixxtest} 
\end{equation}
Transforming the first term on the left hand side as 
\begin{gather*}
-\int\limits_{D_t} \psi_{\tau}\sum\limits_{i=1}^d(a_i\psi_{x_i})_{x_i}\,dx\,d\tau = 
\frac{1}{2}\int\limits_\Omega \sum\limits_{i=1}^d a_i(x,t) \psi_{x_i}^2(x,t) \,dx  - \frac{1}{2} \int\limits_{D_t} \sum\limits_{i=1}^d (a_{i})_\tau \psi_{x_i}^2 \,dx \,d\tau, 
\end{gather*}
and using \eqref{apositive}, from \eqref{psixxtest}, we derive
\begin{gather}
	\frac{a_0}{2}\int\limits_\Omega |D\psi(x,t)|^2 \,dx  +\int\limits_{D_t} z^\ep\Big(\sum\limits_{i=1}^d(a_i\psi_{x_i})_{x_i}\Big)^2 \,dx \,d\tau \nonumber \\
	\leq \int\limits_{D_t} z^\ep\Big(\sum\limits_{i=1}^d (b_i-c_i)\psi_{x_i}\Big)\Big(\sum\limits_{i=1}^d(a_i\psi_{x_i})_{x_i} \Big) \,dx \,d\tau + \int\limits_{D_t} z^\ep\Big((r-\sum\limits_{i=1}^d c_{i,x_i})\psi\Big)\Big(\sum\limits_{i=1}^d(a_i\psi_{x_i})_{x_i} \Big) \,dx \,d\tau \nonumber \\
	+\int\limits_{D_t} \sum\limits_{i=1}^da_i\psi_{x_i}F_{x_i}\,dx\,d\tau + \frac{1}{2}\int\limits_{D_t} \sum\limits_{i=1}^d (a_{i})_\tau \psi_{x_i}^2 \,dx \,d\tau, \label{psixxtest2}
\end{gather}
where $D\psi$ denotes the spatial gradient of $\psi$. Using Cauchy inequality with appropriately chosen small parameter, and \eqref{databounds}, we get the following estimations:
\begin{gather*}
	\int\limits_{D_t} z^\ep\Big(\sum\limits_{i=1}^d (b_i-c_i)\psi_{x_i}\Big)\Big(\sum\limits_{i=1}^d(a_i\psi_{x_i})_{x_i} \Big) \,dx \,d\tau 
    \leq \frac{1}{4}\int\limits_{D_t} z^\ep\Big(\sum\limits_{i=1}^d(a_i\psi_{x_i})_{x_i} \Big)^2 \,dx \,d\tau \\ + 2d\bar{b}^{-1}\Big(\max\limits_{i} \|b_i\|^2_{L_\infty(D)}+\max\limits_{i}\|c_i\|^2_{L_\infty(D)}\Big)t \, \|D\psi\|_{L_{2,\infty}(D_t)}^2,
\end{gather*}
\begin{gather*}
	\int\limits_{D_t} z^\ep\Big((r-\sum\limits_{i=1}^d c_{i,x_i})\psi\Big)\Big(\sum\limits_{i=1}^d(a_i\psi_{x_i})_{x_i} \Big) \,dx \,d\tau 
	\leq \frac{1}{4}\int\limits_{D_t} z^\ep\Big(\sum\limits_{i=1}^d(a_i\psi_{x_i})_{x_i} \Big)^2 \,dx \,d\tau\\ + 2\bar{b}^{-1}\Big(\|r\|^2_{L_\infty(D)}+ d\sum\limits_{i=1}^d\|c_{i,x_i}\|^2_{L_\infty(D)}\Big) \|\psi\|_{L_2(D_t)}^2,  
\end{gather*}
\begin{gather*}
	\int\limits_{D_t} \sum\limits_{i=1}^da_i\psi_{x_i}F_{x_i}\,dx\,d\tau 
    \leq t \|D\psi\|_{L_{2,\infty}(D_t)}^2  + \frac{1}{4}\max\limits_{i}\|a_i\|^2_{L_\infty(D)} \|DF\|_{L_2(D_t)}^2,
\end{gather*}
\begin{gather*}
	\frac{1}{2}\int\limits_{D_t} \sum\limits_{i=1}^d (a_{i})_\tau \psi_{x_i}^2 \,dx \,d\tau \leq  \frac{1}{2}\max\limits_{i}\Big \|\frac{\partial a_i}{\partial t}\Big\|_{L_{\infty,1}(D_t)} \|D\psi\|_{L_{2,\infty}(D_t)}^2. 
\end{gather*}
Plugging these estimates in \eqref{psixxtest2}, absorbing similar terms to the left hand side, and
by taking $\esssup$ with respect to $\tau$ in $0\leq \tau \leq t$ in the first term, we have
\begin{gather}
	\frac{a_0}{2}\|D\psi\|_{L_{2,\infty}(D_t)}^2  +\frac{1}{2}\int\limits_{D_t} z^\ep\Big(\sum\limits_{i=1}^d(a_i\psi_{x_i})_{x_i}\Big)^2 \,dx \,d\tau \leq \nonumber \\
     \Bigg(\Big(2d\bar{b}^{-1}\big(\max\limits_{i} \|b_i\|^2_{L_\infty(D)}+\max\limits_{i}\|c_i\|^2_{L_\infty(D)})+1\Big)t+ \frac{1}{2}\max\limits_{i}\Big \|\frac{\partial a_i}{\partial t}\Big\|_{L_{\infty,1}(D_t)} \Bigg) \|D\psi\|_{L_{2,\infty}(D_t)}^2  \nonumber\\
    + 2\bar{b}^{-1}\Big(\|r\|^2_{L_\infty(D)}+ d\sum\limits_{i=1}^d\|c_{i,x_i}\|^2_{L_\infty(D)}\Big) \|\psi\|_{L_2(D_t)}^2 +  \frac{1}{4}\max\limits_{i}\|a_i\|^2_{L_\infty(D)}  \|DF\|_{L_2(D_t)}^2. \label{psixxtest3}
\end{gather}
By choosing $t>0$ sufficiently small such that 
\begin{equation}\label{tsmall}
   \Big(2d\bar{b}^{-1}\big(\max\limits_{i} \|b_i\|^2_{L_\infty(D)}+\max\limits_{i}\|c_i\|^2_{L_\infty(D)})+1\Big)t+ \frac{1}{2}\max\limits_{i}\Big \|\frac{\partial a_i}{\partial t}\Big\|_{L_{\infty,1}(D_t)} < \frac{a_0}{4},
\end{equation}
and by absorbing the first term on the right hand side, we derive
\begin{gather}
	\frac{a_0}{4}\|D\psi\|_{L_{2,\infty}(D_t)}^2  +\frac{1}{2}\int\limits_{D_t} z^\ep\Big(\sum\limits_{i=1}^d(a_i\psi_{x_i})_{x_i}\Big)^2 \,dx \,d\tau \leq  \nonumber \\
    2\bar{b}^{-1}\Big(\|r\|^2_{L_\infty(D)}+ d\sum\limits_{i=1}^d\|c_{i,x_i}\|^2_{L_\infty(D)}\Big) \|\psi\|_{L_2(D_t)}^2 +  \frac{1}{4}\max\limits_{i}\|a_i\|^2_{L_\infty(D)}  \|DF\|_{L_2(D_t)}^2.  \label{psixxtest4}
\end{gather}
From the maximum principle (e.g. Theorem 2.1, Chapter 1 of \cite{LSU}) it follows that $\psi$ is essentially bounded, and $\|\psi\|_{L_\infty(D)}$ depends on $\bar{b}$ and $L_\infty$-norms of $a_i, (a_i)_{x_i}, b_i, c_i, (c_i)_{x_i}, r$ and $F$. Hence $\|\psi\|_{L_2(D)}$, and therefore the right hand side of \eqref{psixxtest4} is bounded uniformly with respect to $\ep$. If \eqref{tsmall} is not satisfied in the whole time interval $[0,T]$, it can be divided into finitely many intervals that satisfy \eqref{tsmall}, and summing up respective inequalities \eqref{psixxtest4} we arrive at the estimation:
\begin{equation}\label{uniformbound}
	\|D\psi\|_{L_{2,\infty}(D)}^2  +\int\limits_D z^\ep\Big(\sum\limits_{i=1}^d(a_i\psi_{x_i})_{x_i}\Big)^2 \,dx \,d\tau \leq C
\end{equation}
 where $C$ is independent of $\ep$. Furthermore, any constant independent of $\ep$ will be denoted by $C$. Although \eqref{uniformbound} is satisfactory for our purpose, it is worth mentioning that since $\psi_t =F-z^\ep \mathcal{L}^*\psi$,  
from \eqref{databounds},\eqref{zbound},\eqref{uniformbound} and $L_\infty$ bounds of $\psi$ and $F$ it follows that  $\Vert\psi_t\Vert_{L_2(D)}$ is uniformly bounded.
Therefore, complete energy estimate for $\psi$ reads
\begin{equation}\label{uniformboundcomplete}
	\|\psi\|_{L_\infty(D)}+\|\psi_t\|^2_{L_2(D)}+\|D\psi\|_{L_{2,\infty}(D)}^2  +\int\limits_D z^\ep\Big(\sum\limits_{i=1}^d(a_i\psi_{x_i})_{x_i}\Big)^2 \,dx \,d\tau \leq C
\end{equation}
Having \eqref{uniformboundcomplete}, we can estimate $\ep$-term in \eqref{almost} as follows:
\begin{gather*}
    \left| \int\limits_D \hat{B} \ep \mathcal{L}^*\psi\,dx\,dt \right|
    \leq \esssup\limits_{D} |\hat{B}| \Bigg[ \Big( \int\limits_D \frac{\ep^2}{z^\ep} \,dx\,dt  \int\limits_D z^\ep \Big(\sum\limits_{i=1}^d(a_i\psi_{x_i})_{x_i}\Big)^2 \,dx\,dt\Big)^\frac{1}{2}  +C \ep\Bigg] \\
    \leq C \esssup\limits_{D} |\hat{B}| \Big(\ep^{\frac{1}{2}} |D|^{\frac{1}{2}}+\ep\Big ) \to 0 \ \text{as} \ \ep \to 0.
\end{gather*}
Therefore, (\ref{almost}) implies
\[
\int\limits_D \hat{B} F\,dx\,dt = 0.
\]
Since the choice of $F$ is arbitrary, it follows that
\[ B(x,t,v(x,t)) = \tilde B(x,t,\tilde v(x,t)), \ \text{a.e.} \ (x,t)\in D. \] 
Since, $B,\tilde{B} \in \n B$, clearly we have
\[ \{(x,t)\in D: B(x,t,v(x,t)) = \tilde B(x,t,\tilde v(x,t))\} \subset \{(x,t)\in D: v(x,t)=\tilde{v}(x,t)\} \]
and therefore, $v(x,t)=\tilde v(x,t) \ \text{for a.e.}~(x,t) \in D$. Uniqueness is proved. \hfill{$\square$}

The following lemma recalls the criteria for the convergence of the discrete optimal control problems.
\begin{lemma}\label{Vasil}\cite{Vasilev1} The sequence of discrete optimal control problems $\m I_{n}$ approximates the continuous optimal control problem $\m I$ with respect to the functional, i.e. \eqref{approx1} holds, if and only if the following conditions are satisfied:
	\begin{description}
		\item{(i)} For any $f\in\n F^R$, we have $\n Q_{\D}(f)\in\n F_{\D}^R$, and 
		\begin{equation}\label{firstcondition}
			\limsup\limits_{\D\ra0}\Big( \n I_{\D}(\n{Q}_{\D}(f))-\n{J}(f) \Big ) \le 0.
		\end{equation}
		\item{(ii)} For any $[f]_{\D}\in\n F_{\D}^R$, we have $\n P_{\D}([f]_{\D})\in\n F^R$, and
		\begin{equation}\label{secondcondition}
			\limsup\limits_{\D\ra0} \Big ( \n{J}(\n{P}_{\D}([f]_{\D})) -\n{I}_{\D}([f]_{\D})  \Big ) \le 0.
		\end{equation}
		\end{description}
\end{lemma}
Finally, we recall two lemmas proved in \cite{Abdulla7}.
\begin{lemma}\label{PQ}\cite{Abdulla7} The maps $\n P_{\D}$ and $\n Q_{\D}$ satisfy the conditions of Lemma \ref{Vasil}.
\end{lemma}
\begin{lemma}\label{phistrong}\cite{Abdulla7} For any fixed $\ep>0$, there exists $\delta>0$ such that
    \begin{equation}\label{phibound}
    \sum\limits_{\n A}h^d\sum\limits_{i=1}^d|\Phi_{\gamma x_i}|^2\leq(1+\ep)\Vert D\Phi\Vert_{L_2(\Omega)}^2
    \end{equation}
    whenever $h<\delta$.
\end{lemma}

\section{Discrete Energy Estimates}\label{energyestimates}

\begin{theorem}\label{boundedness} (Discrete Maximum Principle) For any $R>0$, $[f]_{\D}\in\n F_{\D}^R$, and $\D$, the discrete state vector $[v([f]_{\D})]_{\D}$ given in Definition \ref{dsvdef} satisfies
\begin{equation}\label{boundedest}
\Vert[v]_{\D}\Vert_{\ell_{\infty}}
\leq e^{\lambda T} \max\left\{\Vert[f]_{\D}\Vert_{\ell_{\infty}}~,~\Vert\Phi\Vert_{L_{\infty}(\Omega)}\right\}
\end{equation}
where 
\begin{equation}\label{lambda}
\lambda = \frac{2}{\bar{b}}\Big(1+ \sum\limits_{i=1}^d\|b_{i,x_i}\|_{L_\infty(D)}+\|r\|_{L_\infty(D)}\Big) 
\end{equation}
\end{theorem}

\emph{Proof.} Fix a discretization $\D=(\tau,h)$ and $[f]_{\D}\in\n F_{\D}^R$. By Lemma \ref{exuniqdsv}, there exists a unique discrete state vector, $[v([f]_{\D})]_{\D}$. We transform it using \eqref{lambda} as
\begin{equation}\label{udef}
u_{\gamma}(k):=v_{\gamma}(k)e^{-\lambda t_k},\quad\forall(\gamma,k)\in\n A(D_{\D}).
\end{equation}
Then by \eqref{mvt}, we get
\begin{align*}
\frac{b_n(v_{\gamma}(k))-b_n(v_{\gamma}(k-1))}{\tau}=\zeta_{\D}^{\gamma,k}e^{\lambda t_{k-1}}u_{\gamma \bar{t}}(k)+\zeta_{\D}^{\gamma,k} \lambda e^{\lambda t^k}u_{\gamma}(k)
\end{align*}
where $t^k\in[t_{k-1},t_k]$ represents the value resulting from the mean value theorem:
\[
e^{\lambda t_k}-e^{\lambda t_{k-1}}=\lambda e^{\lambda t^k}\tau.
\]
Substituting in \eqref{dsvequiv}, we get
\begin{gather*}\label{dsvu}
    \zeta_{\D}^{\gamma,k}\lambda e^{\lambda t^k} u_{\gamma}(k)+\zeta_{\D}^{\gamma,k}u_{\gamma\bar{t}}(k)\,e^{\lambda t_{k-1}}-e^{\lambda t_k}\sum\limits_{i=1}^{d}\Big([(a_{i})_\alpha u_{\gamma x_i}(k)+(b_i)_\alpha u_\gamma(k)]\Big)_{\bar{x}_i}  \\
    +e^{\lambda t_k}\sum\limits_{i=1}^{d} (c_i)_\alpha u_{\gamma x_i}(k) ~+e^{\lambda t_k} r_\alpha u_\gamma(k)=f_{\gamma,k}^{\D}
\end{gather*}
Splitting up the third term and gathering similar terms, we have
\begin{gather}
    \zeta_{\D}^{\gamma,k} \lambda e^{\lambda t^k}u_{\gamma}(k) +\zeta_{\D}^{\gamma,k} e^{\lambda t_{k-1}}u_{\gamma \bar{t}}(k)\,e^{\lambda t_{k-1}} \nonumber \\
    -e^{\lambda t_k}\sum\limits_{i=1}^{d}\frac{1}{h}\Big(((a_{i})_\alpha - h(c_i)_\alpha) u_{\gamma x_i}(k) - ((a_i)_{\alpha - e_i}-h(b_i)_{\alpha-e_i})u_{\gamma \bar{x}_i}(k)\Big) \nonumber\\
    +e^{\lambda t_k}\Big(r-\sum\limits_{i=1}^{d} (b_i)_{\alpha, \bar{x_i}}\Big) u_{\gamma}(k) =f_{\gamma,k}^{\D} \label{dsvu2}
\end{gather}
If $u_{\gamma}(k)\leq 0$ for every $\alpha\in\n A(D_{\D})$, then it is clear that $\max\limits_{\n A(D_{\D})}u_{\gamma}(k)\leq0$. We now suppose that for some $\alpha=(\gamma,k)\in\n A(D_{\D})$, we have $u_{\gamma}(k)>0$. This implies that $\max\limits_{\n A(D_{\D})}u_{\gamma}(k)>0$. Assume that maximum occurs as $\alpha^*=(\gamma^*,k^*)$, i.e.
\[
u_{\gamma^*}(k^*)=\max\limits_{\n A(D_{\D})}u_{\gamma}(k).
\]
Due to (iii) in Definition \ref{dsvdef}, $\alpha^* \notin \n A(S_{\D})$. If $\alpha^*=(\gamma^*,0), \,\gamma^* \in \n A(\Omega'_{\D})$, this would imply
\[
u_{\gamma^*}(k^*)=\max\limits_{\n A(\Omega_{\D})}\Phi_{\gamma}\leq\Vert\Phi\Vert_{L_{\infty}(\Omega)}.
\]
The only other possibility is $\alpha^* \in \n A(D'_{\D})$, i.e. \eqref{dsvu2} is true for $\alpha^*$, and moreover,
\[
u_{\gamma^*\bar{t}}(k^*)\geq0,\qquad u_{\gamma^*x_i}(k^*)\leq0~\forall i,\qquad u_{\gamma^*\bar{x}_{i}}(k^*)\geq0~\forall i
\]
since maximum occurs at $\alpha^*$. By using this properties and \eqref{zetabig} in \eqref{dsvu2}, we have 
\begin{gather}
	e^{\lambda t_{k^*}} \Bigg( \lambda \bar{b}e^{-\lambda(t_{k^*}-t^{k^*})} - \sum\limits_{i=1}^d(b_i)_{\alpha^*,\bar{x}_i} +r_{\alpha^*}\Bigg) u_{\gamma^*}(k^*) \nonumber \\
    -e^{\lambda t_{k^*}}\sum\limits_{i=1}^{d}\frac{1}{h}\Big(((a_{i})_{\alpha^*} - h(c_i)_{\alpha^*}) u_{\gamma^* x_i}(k) - ((a_i)_{\alpha^* - e_i}-h(b_i)_{\alpha^*-e_i})u_{\gamma^* \bar{x}_i}(k)\Big) 
     \leq f_{\gamma^*,k^*}^{\D} \label{dsvu3}
\end{gather}
Assume that $\forall i = 1,...,d$, we have $ (a_{i})_{\alpha^*} - h(c_i)_{\alpha^*}-(a_i)_{\alpha^* - e_i}+h(b_i)_{\alpha^*-e_i} \geq 0$. Then we can rewrite \eqref{dsvu3} as
\begin{gather}
	e^{\lambda t_{k^*}} \Bigg( \lambda \bar{b}e^{-\lambda(t_{k^*}-t^{k^*})} - \sum\limits_{i=1}^d(b_i)_{\alpha^*,\bar{x}_i} +r_\alpha^*\Bigg) u_{\gamma^*}(k^*) \nonumber \\
	-e^{\lambda t_{k^*}}\sum\limits_{i=1}^{d}\frac{1}{h}\Big((a_{i})_{\alpha^*} - h(c_i)_{\alpha^*} - (a_i)_{\alpha^* - e_i}+h(b_i)_{\alpha^*-e_i}\Big)u_{\gamma^* x_i}(k) \nonumber \\
	-e^{\lambda t_{k^*}}\sum\limits_{i=1}^d\Big((a_i)_{\alpha^* - e_i}-h(b_i)_{\alpha^*-e_i}\Big)u_{\gamma^* x_i \bar{x}_i}(k^*)\leq f_{\gamma^*,k^*}^{\D} \label{dsvu4}
\end{gather}
Since $u_{\gamma^* x_i(k^*)} \leq 0$, $u_{\gamma^* x_i \bar{x}_i (k^*)} \leq 0$ and since for small enough $h$, we have $(a_i)_{\alpha^* - e_i}-h(b_i)_{\alpha^*-e_i} \geq 0$, from \eqref{dsvu4} we deduce
\begin{gather}\label{dsvu5}
	e^{\lambda t_{k^*}} \Bigg( \lambda \bar{b}e^{-\lambda(t_{k^*}-t^{k^*})} - \sum\limits_{i=1}^d(b_i)_{\alpha^*,\bar{x}_i} +r_\alpha^*\Bigg) u_{\gamma^*}(k^*) \leq f_{\gamma^*,k^*}^{\D}
\end{gather}
If for some $i = \overline{1,d}$ we have that $ (a_{i})_{\alpha^*} - h(c_i)_{\alpha^*}-(a_i)_{\alpha^* - e_i}+h(b_i)_{\alpha^*-e_i} < 0$, then we rewrite that specific term in  \eqref{dsvu3} as 
\begin{gather*}
	\frac{1}{h}\Bigg(((a_{i})_{\alpha^*} - h(c_i)_{\alpha^*}) u_{\gamma^* x_i}(k) - ((a_i)_{\alpha^* - e_i}-h(b_i)_{\alpha^*-e_i})u_{\gamma^* \bar{x}_i}(k)\Bigg) \\
	= \frac{1}{h}\Big((a_{i})_{\alpha^*} - h(c_i)_{\alpha^*} - (a_i)_{\alpha^* - e_i}+h(b_i)_{\alpha^*-e_i}\Big)u_{\gamma^* \bar{x}_i}(k) +((a_{i})_{\alpha^*} - h(c_i)_{\alpha^*})u_{\gamma^* x_i \bar{x}_i}(k^*)
\end{gather*}
Since $u_{\gamma^* \bar{x}_i(k^*)} \geq 0$, $u_{\gamma^* x_i \bar{x}_i (k^*)} \leq 0$ and since for small enough $h$, we have $(a_i)_{\alpha^*}-h(c_i)_{\alpha^*-e_i} \geq 0$, we get the same estimate as \eqref{dsvu5}. By choosing $\tau$ sufficiently small such that $e^{-\lambda(t_{k}-t^{k})}>\frac12\quad,\forall k$, due to \eqref{lambda} we get 
\[
u_{\gamma^*}(k^*)\leq f_{\alpha^*}^{\D}e^{-\lambda t^k}\leq\Vert[f]_{\D}\Vert_{\ell_{\infty}}.
\]
These estimations result in
\begin{gather*}
\max\limits_{\n A(D_{\D})}v_{\gamma}(k)\leq e^{\lambda T}\max\left\{\Vert[f]_{\D}\Vert_{\ell_{\infty}}~,~\Vert\Phi\Vert_{L_{\infty}(\Omega)}\right\}.
\end{gather*}
We can similarly derive a uniform lower bound
\begin{gather*}
\min\limits_{\n A(D_{\D})}v_{\gamma}(k)\geq e^{\lambda T}\min\left\{-\Vert[f]_{\D}\Vert_{\ell_{\infty}}~,~-\Vert\Phi\Vert_{L_{\infty}(\Omega)}\right\},
\end{gather*}
giving \eqref{boundedest}, and thus proving the theorem.\hfill{$\square$}\\

\begin{theorem}\label{energy}(Discrete $W_2^{1,1}$ Energy Estimate) For any $R>0$, $\D$ and $[f]_{\D}\in\n F_{\D}^R$, the discrete state vector $[v([f]_{\D})]_{\D}$ satisfies
	\begin{gather}
	\sum\limits_{k=1}^n\tau\sum\limits_{\n A}h^d (v_{\gamma\bar t}(k))^2+\max\limits_{1\leq k\leq n}\sum\limits_{\n A}h^d\sum\limits_{i=1}^d (v_{\gamma x_i}(k))^2+\nonumber\\[4mm]+\sum\limits_{k=1}^n\tau^2\sum\limits_{\n A}h^d\sum\limits_{i=1}^d \big(v_{\gamma x_i\bar t}(k)\big)^2\leq C~\Bigg(\Vert \Phi \Vert_{L_\infty(\Omega)}^2+ \Vert D\Phi \Vert_{L_2(\Omega)}^2+\Vert f^{\D}\Vert_{L_\infty(D)}^2\Bigg)\label{energyest}
	\end{gather}
	where $C$ is a constant independent of $\D, R$.
\end{theorem}
\emph{Proof.} Choose $\eta_\gamma=2\tau v_{\gamma \bar{t}}(k)$ in \eqref{dsveq} with $k=1,...,n$. Using \eqref{mvt} and the identity
\begin{gather*}
2\tau (a_i)_\alpha v_{\gamma x_i}(k)\big(v_{\gamma}(k)_{\bar t}\big)_{x_i}\\
=(a_i)_\alpha (v_{\gamma x_i}(k))^2-(a_i)_{\alpha-e_k}(v_{\gamma x_i}(k-1))^2+\tau^2(a_i)_\alpha \big(v_{\gamma x_i\bar t}(k)\big)^2-\tau (a_i)_{\alpha\bar{t}}(v_{\gamma x_i}(k-1))^2
\end{gather*}
we have
\begin{gather}
    \sum\limits_{\n A}h^d\Bigg[2\tau\zeta_{\D}^{\gamma,k}(v_{\gamma\bar t}(k))^2+\sum\limits_{i=1}^d (a_i)_\alpha (v_{\gamma x_i}(k))^2-\sum\limits_{i=1}^d(a_i)_{\alpha-e_k}(v_{\gamma x_i}(k-1))^2 \nonumber \\
    +\sum\limits_{i=1}^d\tau^2(a_i)_\alpha \big(v_{\gamma x_i\bar t}(k)\big)^2-\sum\limits_{i=1}^d\tau (a_i)_{\alpha\bar{t}}(v_{\gamma x_i}(k-1))^2\Big)+\sum\limits_{i=1}^d2\tau(b_i)_\alpha v_\gamma(k)v_{\gamma x_i\bar t}(k) \nonumber \\
    \sum\limits_{i=1}^d2\tau(c_i)_\alpha v_{\gamma x_i}(k)v_{\gamma \bar t}(k)+2\tau \,r_\alpha v_{\gamma}(k)v_{\gamma \bar t}(k)-2\tau f_{\gamma,k}^{\D}v_{\gamma\bar t}(k)\Bigg]=0. \label{eq1}
\end{gather}
Since $v_{\gamma}(k)=0$ for $\gamma\in\n A(\partial\Omega_{\D})$ through summation by parts we deduce
\begin{gather}
	\sum\limits_{\n A}h^d \sum\limits_{i=1}^d2\tau(b_i)_\alpha v_\gamma(k)v_{\gamma x_i\bar t}(k) = -\sum\limits_{\n A}h^d \sum\limits_{i=1}^d2\tau((b_i)_\alpha v_\gamma(k))_{\bar{x}_i} v_{\gamma\bar t}(k)
	\nonumber\\ = -\sum\limits_{\n A}h^d \sum\limits_{i=1}^d2\tau(b_i)_{\alpha \bar{x}_i} v_\gamma(k)v_{\gamma\bar t}(k)- \sum\limits_{\n A}h^d \sum\limits_{i=1}^d2\tau(b_i)_{\alpha-e_i} v_{\gamma-e_i,x_i}(k)v_{\gamma\bar t}(k). \label{sumparts1}
\end{gather}
Due to (\ref{zetabig}) and \eqref{sumparts1}, from \eqref{eq1} it follows
\begin{gather}
	\sum\limits_{\n A}h^d\Big[2\tau\bar{b}(v_{\gamma\bar t}(k))^2+\sum\limits_{i=1}^d (a_i)_\alpha (v_{\gamma x_i}(k))^2-\sum\limits_{i=1}^d(a_i)_{\alpha-e_k}(v_{\gamma x_i}(k-1))^2 \nonumber \\
	+\sum\limits_{i=1}^d\tau^2(a_i)_\alpha \big(v_{\gamma x_i\bar t}(k)\big)^2 \Big]\leq \sum\limits_{\n A}h^d\Big[ \sum\limits_{i=1}^d\tau (a_i)_{\alpha\bar{t}}(v_{\gamma x_i}(k-1))^2 + \sum\limits_{i=1}^d2\tau(b_i)_{\alpha \bar{x}_i} v_\gamma(k)v_{\gamma\bar t}(k)\nonumber \\
  +  \sum\limits_{i=1}^d2\tau(b_i)_{\alpha-e_i} v_{\gamma-e_i,x_i}(k)v_{\gamma\bar t}(k) -\sum\limits_{i=1}^d2\tau(c_i)_\alpha v_{\gamma x_i}(k)v_{\gamma \bar t}(k) \nonumber \\
    -2\tau(d)_\alpha v_{\gamma}(k)v_{\gamma \bar t}(k)+2\tau f_{\gamma,k}^{\D}v_{\gamma\bar t}(k)\Big] \label{eq2}
\end{gather}
Applying Cauchy inequality with appropriately chosen small parameter we estimate various terms in \eqref{eq2} as follows:
\begin{gather*}
	 \sum\limits_{i=1}^d2\tau(b_i)_{\alpha-e_i} v_{\gamma-e_i,x_i}(k)v_{\gamma\bar t}(k) 
	\leq \tau \frac{\bar{b}}{5} (v_{\gamma\bar t}(k))^2 +  \frac{5d}{\bar{b}}\max\limits_{i}\|b_i\|^2_{L_\infty(D)} \tau\sum\limits_{i=1}^d(v_{\gamma-e_i,x_i}(k))^2, 
	\end{gather*}
	\begin{gather*}
	 \sum\limits_{i=1}^d2\tau(b_i)_{\alpha \bar{x}_i} v_\gamma(k)v_{\gamma\bar t}(k) 
	\leq  \tau \frac{\bar{b}}{5} (v_{\gamma\bar t}(k))^2 + \frac{5d^2}{\bar{b}} \max\limits_{i}\|(b_i)_{x_i}\|^2_{L_\infty(D)}\tau (v_\gamma(k))^2, 
	\end{gather*}
	\begin{gather*}
	-\sum\limits_{i=1}^d2\tau(c_i)_\alpha v_{\gamma x_i}(k)v_{\gamma \bar t}(k) 
	\leq  \tau \frac{\bar{b}}{5} (v_{\gamma\bar t}(k))^2 + \frac{5d}{\bar{b}}\max\limits_{i}\|(c_i)\|^2_{L_\infty(D)}\tau\sum\limits_{i=1}^d (v_{\gamma x_i}(k))^2,
	\end{gather*}
	\begin{gather*}
	-2\tau \, r_\alpha v_{\gamma}(k)v_{\gamma \bar t}(k) \leq \frac{\bar{b}}{5}(v_{\gamma \bar t}(k))^2 + \frac{5}{\bar{b}}\|r\|^2_{L_\infty(D)}(v_\gamma(k))^2,
	\end{gather*}
	\begin{gather*}
	2\tau f_{\gamma,k}^{\D}v_{\gamma\bar t}(k) \leq \frac{\bar{b}}{5}(v_{\gamma \bar t}(k))^2 + \frac{5}{\bar{b}}\tau (f_{\gamma,k}^{\D})^2
\end{gather*}
Implementing these estimates in \eqref{eq2}, and absorbing similar terms into the left hand side, we derive 
\begin{gather}
	\sum\limits_{\n A}h^d\Bigg[\tau\bar{b}(v_{\gamma\bar t}(k))^2+\sum\limits_{i=1}^d (a_i)_\alpha (v_{\gamma x_i}(k))^2 - \sum\limits_{i=1}^d (a_i)_{\alpha-e_k} (v_{\gamma x_i}(k-1))^2 \nonumber \\
	+\sum\limits_{i=1}^d\tau^2(a_i)_\alpha \big(v_{\gamma x_i\bar t}(k)\big)^2 \Bigg]
	\leq \sum\limits_{\n A}h^d\Bigg[ \sum\limits_{i=1}^d\tau (a_i)_{\alpha\bar{t}}(v_{\gamma x_i}(k-1))^2 \nonumber\\
	+\frac{5d}{\bar{b}}(\max\limits_i \|b_i\|^2_{L_\infty(D)}+\max\limits_i \|c_i\|^2_{L_\infty(D)} )\tau \sum\limits_{i=1}^d (v_{\gamma x_i}(k))^2 \nonumber \\
	+ \frac{5d^2}{\bar{b}}\max\limits_i \|(b_i)_{x_{i}}\|^2_{L_\infty(D)}+\frac{5}{\bar{b}}\|r\|^2_{L_\infty(D)}\tau(v_\gamma(k))^2+\frac{5}{\bar{b}}\tau (f_{\gamma,k}^{\D})^2 \Bigg] \label{eq3}.
\end{gather}
Pursuing summation over all $k=\overline{1,q}$, where $q \leq n$, and using \eqref{apositive}, from \eqref{eq3} we have
\begin{gather}
	\sum\limits_{k=1}^q \tau \sum\limits_{\n A}h^d \bar{b}(v_{\gamma\bar t}(k))^2+ a_0\sum\limits_{\n A}h^d\sum\limits_{i=1}^d  (v_{\gamma x_i}(q))^2 \nonumber \\ 
	+a_0 \sum\limits_{k=1}^q \tau^2 \sum\limits_{\n A}h^d\sum\limits_{i=1}^d\big(v_{\gamma x_i\bar t}(k)\big)^2  \leq \sum\limits_{k=1}^q\tau\sum\limits_{\n A}h^d\sum\limits_{i=1}^d (a_i)_{\alpha\bar{t}}(v_{\gamma x_i}(k-1))^2 \nonumber \\
	+\frac{5d}{\bar{b}}\Big(\max\limits_i \|b_i\|^2_{L_\infty(D)}+\max\limits_i \|c_i\|^2_{L_\infty(D)}\Big ) \sum\limits_{k=1}^q \tau\sum\limits_{\n A}h^d \sum\limits_{i=1}^d (v_{\gamma x_i}(k))^2 \nonumber \\
	+ \frac{5}{\bar{b}}\Big(d^2\max\limits_i \|(b_i)_{x_{i}}\|^2_{L_\infty(D)}+\|r\|^2_{L_\infty(D)}\Big )\sum\limits_{k=1}^q \tau \sum\limits_{\n A}h^d(v_\gamma(k))^2\nonumber\\+\frac{5}{\bar{b}}\sum\limits_{k=1}^q\tau\sum\limits_{\n A}h^d(f_{\gamma,k}^{\D})^2 
	+ \sum\limits_{\n A}h^d\sum\limits_{i=1}^d(a_i)_{(\gamma,0)}(v_{\gamma x_i}(0))^2 \label{eq4}
\end{gather}
We estimate the first term on the right hand side as follows:
\begin{gather}
	\sum\limits_{k=1}^q\tau \sum\limits_{\n A}h^d \sum\limits_{i=1}^d(a_i)_{\alpha\bar{t}}(v_{\gamma x_i}(k-1))^2 
	= \sum\limits_{k=0}^{q-1}\tau \sum\limits_{\n A}h^d \sum\limits_{i=1}^d(a_i)_{\alpha t}(v_{\gamma x_i}(k))^2 \nonumber \\
    = \sum\limits_{k=1}^{q-1} \sum\limits_{\n A}\sum\limits_{i=1}^d \Bigg( \frac{1}{\tau} \int\limits_{R_{\D}^\gamma} \int\limits_{t_{k}}^{t_{k+1}} \int\limits_{t-\tau}^{t} \frac{\partial a(x,\xi)}{\partial \xi} \,d\xi \, dt \, dx\Bigg) (v_{\gamma x_i}(k))^2 \nonumber \\
    +\sum\limits_{\n A} \sum\limits_{i=1}^d \Bigg(\frac{1}{\tau} \int\limits_{R_{\D}^\gamma} \int\limits_{0}^{\tau} \int\limits_{0}^{t} \frac{\partial a(x,\xi)}{\partial \xi} \,d\xi \, dt \, dx \Bigg)(v_{\gamma x_i}(0))^2 \leq \nonumber \\
    \leq 2 \max\limits_{1 \leq i \leq d} \Big\|\frac{\partial a_i}{\partial t}\Big\|_{L_{\infty,1}(D_{t_q})}  \max\limits_{1\leq k \leq q} \sum\limits_{\n A}h^d \sum\limits_{i=1}^d(v_{\gamma x_i}(k))^2 \nonumber \\
    +\max\limits_{1 \leq i \leq d} \Big\|\frac{\partial a_i}{\partial t}\Big\|_{L_{\infty,1}(D_{\tau})} \sum\limits_{\n A}h^d \sum\limits_{i=1}^d (\Phi_{\gamma x_i})^2. \label{314}
\end{gather}
We also have
\begin{gather}
	\sum\limits_{k=1}^q \tau\sum\limits_{\n A}h^d \sum\limits_{i=1}^d (v_{\gamma x_i}(k))^2 \leq t_q \max\limits_{1\leq k \leq q} \sum\limits_{\n A}h^d  \sum\limits_{i=1}^d (v_{\gamma x_i}(k))^2.\label{315}
\end{gather}
Applying \eqref{314} and \eqref{315} in \eqref{eq4}, and by noting that the index $q$ in the second term on the left hand side of \eqref{eq4} can be replaced with any $1\leq k \leq q$, we derive
\begin{gather}
	\bar{b} \sum\limits_{k=1}^q \tau \sum\limits_{\n A}h^d (v_{\gamma\bar t}(k))^2 + a_0 \max\limits_{1\leq k \leq q} \sum\limits_{\n A}h^d \sum\limits_{i=1}^d (v_{\gamma x_i}(k))^2 \nonumber \\ 
	+a_0 \sum\limits_{k=1}^q \tau^2 \sum\limits_{\n A}h^d\sum\limits_{i=1}^d\big(v_{\gamma x_i\bar t}(k)\big)^2  \leq 2 \max\limits_{1 \leq i \leq d} \Big\|\frac{\partial a_i}{\partial t}\Big\|_{L_{\infty,1}(D_{t_q})}\max\limits_{1\leq k \leq q} \sum\limits_{\n A}h^d \sum\limits_{i=1}^d(v_{\gamma x_i}(k))^2 \nonumber \\
	+\frac{5d}{\bar{b}}\Big(\max\limits_i \|b_i\|^2_{L_\infty(D)}+\max\limits_i \|c_i\|^2_{L_\infty(D)}\Big ) t_q \max\limits_{1\leq k \leq q} \sum\limits_{\n A}h^d \sum\limits_{i=1}^d (v_{\gamma x_i}(k))^2 \nonumber \\
	+ \frac{5}{\bar{b}}\Big(d^2\max\limits_i \|(b_i)_{x_{i}}\|^2_{L_\infty(D)}+\|r\|^2_{L_\infty(D)}\Big)T |\Omega| \ \|[v]_\Delta\|^2_{\ell_\infty}+\frac{5}{\bar{b}}\sum\limits_{k=1}^q\tau\sum\limits_{\n A}h^d(f_{\gamma,k}^{\D})^2 \nonumber \\
	+ \Big(\max\limits_{1\leq i \leq d}\|a_i\|_{L_\infty(D)}+\max\limits_{1 \leq i \leq d} \Big\|\frac{\partial a_i}{\partial t}\Big\|_{L_{\infty,1}(D)}\Big)\sum\limits_{\n A}h^d\sum\limits_{i=1}^d(\Phi_{\gamma x_i})^2 \label{eq6}
\end{gather}
Note that by the Cauchy-Schwartz inequality, we have
\begin{equation}\label{cbsf}
\sum\limits_{k=1}^q\tau\sum\limits_{\n A}h^d(f_{\gamma,k}^{\D})^2\leq\int\limits_{D_{\D}}(f^{\D})^2\,dx\,dt\leq\Vert{f^{\D}}\Vert_{L_2(D)}^2.
\end{equation}
If the length of the time interval $T$ is small enough to guarantee
\begin{equation}\label{smallsmall}
2 \max\limits_{1 \leq i \leq d} \Big\|\frac{\partial a_i}{\partial t}\Big\|_{L_{\infty,1}(D)} + \frac{5d}{\bar{b}}\Big(\max\limits_i \|b_i\|^2_{L_\infty(D)}+\max\limits_i \|c_i\|^2_{L_\infty(D)}\Big ) T \leq \frac{a_0}{2},
\end{equation}
then by choosing $q=n$, and by absorbing first two terms on the right hand side of \eqref{eq6} into the second term on the left hand side, and by using \eqref{boundedest},\eqref{phibound},\eqref{cbsf}, from \eqref{eq6}, the energy estimate \eqref{energyest} follows. If \eqref{smallsmall} is not satisfied, then we can partition $[0,T]$ into finitely many subsegments which obey \eqref{smallsmall}, pursue the energy estimation in each subsegment as before, and through summation achieve the same for \eqref{eq6} in general. Theorem is proved. \hfill{$\square$}

Theorems~\ref{boundedness} and~\ref{energy} imply the following corollary:
\begin{corollary}\label{interp} Let $\{[f]_{\D}\}$ be a sequence of discrete control vectors such that there exists $R>0$ for which $[f]_{\D}\in\n F_{\D}^R$ for each $\D$. The following statements hold:
    \begin{enumerate}[(a)]
        \item The sequences $\{\tilde V_{\D}\},\{V_{\D}\}, \{V_{\D}'\}$ are uniformly bounded in $L_{\infty}(D)$.
        
        \item For each $i\in\{1,\ldots,d\}$, the sequences $\{\tilde V_{\D}^i\}, \{\partial V_{\D}/\partial x_i\},\{\partial V_{\D}'/\partial x_i\}$ are uniformly bounded in $L_2(D)$. Moreover, the sequence $\{\partial V_{\D}'/\partial t\}$ is uniformly bounded in $L_2(D)$.
        
        \item The sequence $\{V_{\D}-V_{\D}'\}$ converges strongly to $0$ in $L_2(D)$ as $\tau\ra0$.
        
        \item {For each $k=1,\ldots,n$, the sequence $\{V_{\D}^k-\tilde V_{\D}(\cdot,t_k)\}$ converges strongly to $0$ in $L_2(\Omega)$ as $h\ra0$.} Furthermore, the sequence $\{\tilde V_{\D}-V_{\D}\}$ converges strongly to $0$ in $L_2(D)$ as $h\ra0$.
        
        \item For each $i\in\{1,2,\ldots,d\}$, the sequence $\{\partial V_{\D}/\partial x_i-\partial V_{\D}'/\partial x_i\}$ converges strongly to $0$ in $L_2(D)$ as $\tau\ra0$.
        
        \item For each $i\in\{1,2,\ldots,d\}$, the sequence $\{\tilde V_{\D}^i-\partial V_{\D}/\partial x_i\}$ converges weakly to $0$ in $L_2(D)$ as $\D\ra0$.
    \end{enumerate}
\end{corollary}
Having estimates \eqref{boundedest},\eqref{energyest}, the proof of the corollary coincides with the proof of corresponding result of \cite{Abdulla7} (Theorem 14, pp. 22-30).

\section{Proofs of the Main Results}\label{approximationtheorem}
The key to complete the proof of main results is the following approximation theorem.
\begin{theorem}\label{approx} Let $R>0$ is fixed, and for the sequence of discrete control vectors $[f]_{\D}\in\n F_{\D}^R$, corresponding sequence of interpolations $\{\n P_{\D}([f]_{\D})\}$ converges weakly to $f$ in $L_2(D)$. Then the sequence of multilinear interpolations $\{V_{\D}'\}$ of associated discrete state vectors converges weakly in $W_2^{1,1}(D)$ to weak solution $v=v(x,t;f)\in\overset{\circ}{W}{}_2^{1,1}(D)\cap L_{\infty}(D)$ of the singular PDE problem \eqref{PDE}-\eqref{vgfuture}.
\end{theorem}
\emph{Proof.} From Theorems~\ref{boundedness},~\ref{energy} and Corollary~\ref{interp} it follows that $\{V_{\D}'\}$ is a uniformly bounded sequence in $W_2^{1,1}(D)\cap L_{\infty}(D)$, and hence it is weakly precompact in $W_2^{1,1}(D)$. Let $v$ be its weak limit point. By the Rellich-Kondrachev compact embedding \cite{Nikolski}, there is a subsequence that converges strongly in $L_2(D)$, and hence further subsequence can be chosen which converges pointwise almost everywhere on $D$. Since $\{V_{\D}'\}$ is a uniformly bounded in $L_\infty(D)$, and subspace $\overset{\circ}{W}{}_2^{1,1}(D)$ is closed in the weak topology of $W_2^{1,1}(D)$, it follows $v\in\overset{\circ}{W}{}_2^{1,1}(D)\cap L_{\infty}(D)$. Next, we prove that $v$ is a weak solution of the singular PDE problem \eqref{PDE}-\eqref{vgfuture}.

Without loss of generality assume that the whole sequence $\{V_{\D}'\}$ converges to $v$, weakly in $W_2^{1,1}(D)$ and pointwise a.e. on $D$. Let $\psi\in\overset{\bullet}{\m C}{}^1(D)$ is a continuously differentiable function on $\bar{D}$, whose support is positive distance away from $S$ and $\Omega\times\{t=T\}$. Due to construction of $D_\D$, there exists a discretization $\D^*$, such that $\overline{\text{supp }\psi}\subset D_{\D}$ for all $\D\leq\D^*$. For $\D\leq\D^*$. We define a discrete vector
\[ [\psi]_{\D}=\{ \psi_{\gamma}^k: \ \psi_{\gamma}^k=\psi(x_\gamma,t_k), \   \alpha=(\gamma,k) \in \n A(D_{\D})\} \]
Note that $ \psi_{\gamma}^n=0$, for all $\gamma\in \n A(\Omega_{\D})$. 
Plugging $\eta_{\gamma}:=\tau\psi_{\gamma}^k$ into \eqref{dsveq}, 
and pursuing summation over $k=\overline{1,n}$, we get
\begin{gather}
\sum\limits_{k=1}^n \tau\sum\limits_{\n A}h^d\Bigg[\big(b_n(v_{\gamma}(k))\big)_{\bar t}\psi_{\gamma}^k+\sum\limits_{i=1}^d \Big((a_{i})_{\alpha} v_{\gamma x_i}(k)+(b_i)_{\alpha} v_\gamma(k)\Big)\psi_{\gamma x_i}^k \nonumber \\
+\sum\limits_{i=1}^d (c_i)_\alpha v_{\gamma x_i}(k) \psi_{\gamma}^k+r_\alpha v_\gamma (k) \psi_{\gamma}^k -f_{(\gamma,k)}^{\D}\psi_{\gamma}^k\Bigg]=0.\label{dsvsum}
\end{gather}
Since
\begin{gather}
    \sum\limits_{k=1}^n\tau\sum\limits_{\n A}h^d\big(b_n(v_{\gamma}(k))\big)_{\bar t}\psi_{\gamma}^k=-\sum\limits_{k=1}^{n-1}\tau\sum\limits_{\n A}h^db_n(v_{\gamma}(k))\psi_{\gamma t}^k-\sum\limits_{\n A}h^db_n(\Phi_{\gamma})\psi_{\gamma}^1\label{sumparts},
\end{gather}
from \eqref{dsvsum}) we have
\begin{gather}
	-\sum\limits_{k=1}^{n-1} \tau\sum\limits_{\n A}h^d b_n(v_{\gamma}(k))\psi_{\gamma t}^k+\sum\limits_{k=1}^n \tau\sum\limits_{\n A}h^d\Bigg[\sum\limits_{i=1}^d \Big((a_{i})_{\alpha} v_{\gamma x_i}(k)+(b_i)_{\alpha} v_\gamma(k)\Big)\psi_{\gamma x_i}^k \nonumber \\
	+\sum\limits_{i=1}^d (c_i)_\alpha v_{\gamma x_i}(k) \psi_{\gamma}^k+r_\alpha v_\gamma (k) \psi_{\gamma}^k -f_{(\gamma,k)}^{\D}\psi_{\gamma}^k\Bigg]-\sum\limits_{\n A}h^db_n(\Phi_{\gamma})\psi_{\gamma}^1=0. \label{dsvsum1}
\end{gather}
We define the following interpolations 
\begin{gather}
\overline\Phi_{\D}\Big|_{R_{\D}^{\gamma}}=\Phi_{\gamma},~~\gamma\in\n A,\qquad\overline\Phi_{\D}\equiv0\text{ elsewhere on } \Omega,\nonumber\\[4mm]
\overline\psi_{\D}\Big|_{C_{\D}^{\alpha}}=\psi_{\gamma}^k,~~\alpha\in\n A(\n C_{\D}^D),\qquad\overline\psi_{\D}\equiv0\text{ elsewhere on }D,\nonumber\\[4mm]
\overline\psi_{\D}^t\Big|_{C_{\D}^{\alpha}}=\psi_{\gamma t}^k,~~\alpha\in\n A(\n C_{\D}^D\backslash\n R_{\D}^{\gamma,n}),\qquad\overline\psi_{\D}^t\equiv0\text{ elsewhere on }D,\nonumber\\[4mm]
\overline\psi_{\D}^i\Big|_{C_{\D}^{\alpha}}=\psi_{\gamma x_i}^k,~~\alpha\in\n A(\n C_{\D}^D),~~k=1,\ldots,n,\qquad\overline\psi_{\D}^i\equiv0\text{ elsewhere on }D.\nonumber
\end{gather}
and rewrite \eqref{dsvsum1} in integral form:
\begin{gather}
 -\sum\limits_{k=1}^{n-1}~\int\limits_{t_{k-1}}^{t_k}\sum\limits_{\n A}\int\limits_{R_{\D}^{\gamma}}b_n(\tilde V_{\D})\overline\psi_{\D}^t\,dx\,dt +  \sum\limits_{k=1}^{n}~\int\limits_{t_{k-1}}^{t_k}\sum\limits_{\n A}\int\limits_{R_{\D}^{\gamma}}\Big[\sum\limits_{i=1}^d \Big(a_i(x,t)\tilde{V}_{\D}^i+b_i(x,t)\tilde{V}_{\D}\Big)\overline\psi_{\D}^i \nonumber \\
    +\sum\limits_{i=1}^d c_i(x,t)\tilde{V}_{\D}^i \overline\psi_{\D}+r(x,t) \tilde{V}_{\D} \overline\psi_{\D} -f^{\D}\overline\psi_{\D}\Big]\,dx\,dt -\sum\limits_{\n A}\int\limits_{R_{\D}^{\gamma}}b_n(\overline\Phi_{\D})\overline\psi_{\D}(x,\tau)\,dx=0,\label{int1}
\end{gather}
which then implies
\begin{gather}
    \int\limits_{D}\Big[-b_n(\tilde V_{\D})\overline\psi_{\D}^t+\sum\limits_{i=1}^d \Big(a_i(x,t)\tilde{V}_{\D}^i+b_i(x,t)\tilde{V}_{\D}\Big)\overline\psi_{\D}^i \nonumber \\
    +\sum\limits_{i=1}^d c_i(x,t)\tilde{V}_{\D}^i \overline\psi_{\D}+r(x,t) \tilde{V}_{\D} \overline\psi_{\D} -f^{\D}\overline\psi_{\D}\Big]\,dx\,dt -\int\limits_{\Omega}b_n(\overline\Phi_{\D})\overline\psi_{\D}(x,\tau)\,dx=0,\label{int1}
\end{gather}
due to $\overline\psi_{\D}^t\equiv0$ on $\Omega\times(T-\tau,T]$. We transform \eqref{int1} as follows:
\begin{gather}
\int\limits_{D}\Big[-b_n(\tilde V_{\D})\frac{\partial \psi}{\partial t}+\sum\limits_{i=1}^d \Big(a_i(x,t)\tilde{V}_{\D}^i+b_i(x,t)\tilde{V}_{\D}\Big)\frac{\partial \psi}{\partial x_i}\nonumber \\
+\sum\limits_{i=1}^d c_i(x,t)\tilde{V}_{\D}^i \psi+r(x,t) \tilde{V}_{\D} \psi -f^{\D}\psi\Big]\,dx\,dt -\int\limits_{\Omega}b_n(\overline\Phi_{\D})\psi(x,0)\,dx + I=0, \label{int2}
\end{gather}
where
\begin{gather}
I=\int\limits_{D}\Big[-b_n(\tilde V_{\D})\left(\overline\psi_{\D}^t-\frac{\partial\psi}{\partial t}\right)+\sum\limits_{i=1}^d\Big(a_i(x,t)\tilde V_{\D}^i+b_i(x,t) \tilde{V}_{\D} \Big)\left(\overline\psi_{\D}^i- \frac{\partial\psi}{\partial x_i}\right) \nonumber \\
+\sum\limits_{i=1}^d c_i(x,t)\tilde{V}_{\D}^i \Big(\overline\psi_{\D}-\psi\Big)+r(x,t) \tilde{V}_{\D} \big(\overline\psi_{\D}-\psi\big)-f^{\D}\big(\overline\psi_{\D}-\psi\big)\Big]\,dx\,dt \nonumber \\
-\int\limits_{\Omega}b_n(\overline\Phi_{\D})\Big(\overline\psi_{\D}(x,\tau)-\psi(x,0)\Big)\,dx.\label{I}
\end{gather}
Since sequences $b_n(\tilde V_{\D})$ and $b_n(\overline\Phi_{\D})$ are uniformly bounded, and the sequences $\overline\psi_{\D},\overline\psi_{\D}^t,\overline\psi_{\D}^i$ converge uniformly on $\overline D$ to the functions $\psi,\partial\psi/\partial t,\partial\psi/\partial x_i$ respectively as $\D\ra0$, it easily follows that $I\ra0$ as $\D\ra0$. 
In \cite{Abdulla7}, it is proved that $b_n(\tilde V_{\D})$, and $b_n(\overline\Phi_{\D})$ are weakly convergent sequences in $L_2(D)$ and $L_2(\Omega)$ respectively, and their weak limits are functions of type $\n B$. Precisely, it is proved that
\begin{gather}
b_n(\tilde V_{\D}) \rightharpoonup \tilde{b}(x,t) \ \text{in} \ L_2(D); \ \ \tilde{b}(x,t)=B(x,t,v(x,t)), \ \text{a.e. in} \ D, \label{weaklimclassB}\\
b_n(\overline\Phi_{\D}) \rightharpoonup \tilde{b}_0(x) \ \text{in} \ L_2(\Omega); \ \ \tilde{b}_0(x)=B_0(x,\Phi(x)), \ \text{a.e. in} \ \Omega, \label{weaklimclassB0}
\end{gather}
where $B$ and $B_0$ are some functions of class $\n B$. Passing to limit as $\D \ra 0$, from \eqref{int2},\eqref{weaklimclassB},\eqref{weaklimclassB0} it follows that
\begin{gather}
	\int\limits_0^T\int\limits_{\Omega}\Big[-\tilde{b}(x,t)\frac{\partial \psi}{\partial t}+\sum\limits_{i=1}^d \Big(a_i(x,t)\frac{\partial v}{\partial x_i}+b_i(x,t)v\Big)\frac{\partial \psi}{\partial x_i}\nonumber \\
	+\sum\limits_{i=1}^d c_i(x,t)\frac{\partial v}{\partial x_i} \psi+r(x,t) v \psi -f\psi\Big]\,dx\,dt -\int\limits_{\Omega}b_0(x)\psi(x,0)\,dx =0.
\end{gather}
Since $\overset{\bullet}{\m C}{}^1(D)$ is dense in the set of admissible test functions $\psi$, and by using \eqref{weaklimclassB},\eqref{weaklimclassB0} again, it follows that $v$ is a weak solution of the singular PDE problem \eqref{PDE}-\eqref{vgfuture}. \hfill{$\square$}

Theorem~\ref{approx} and Proposition~\ref{some} together with energy estimates of Section~\ref{energyestimates} imply the general existence, uniqueness and stability result for the singular PDE problem \eqref{PDE}-\eqref{vgfuture}, when the data satisfy assumptions formulated in Section~\ref{mainresults} and $f\in L_\infty(D)$.
\begin{corollary}\label{vestimates} There exists a unique weak solution $v\in\overset{\circ}{W}{}_2^{1,1}(D)\cap L_{\infty}(D)$ of the singular PDE problem \eqref{PDE}-\eqref{vgfuture} and the following estimates are satisfied:
    \begin{equation}\label{boundv}
    \Vert v\Vert_{L_{\infty}(D)}\leq e^{\lambda T}\max\left\{\Vert f\Vert_{L_{\infty}(D)}~,~\Vert\Phi\Vert_{L_{\infty}(\Omega)}\right\},
    \end{equation}
    \begin{equation}\label{energyv}
    \Vert D_xv\Vert_{L_2(D)}^2+\Vert v_t\Vert_{L_2(D)}^2\leq C\left[~\Vert f\Vert_{L_\infty(D)}^2+\Vert\Phi\Vert_{L_\infty(\Omega)}^2+\Vert D\Phi\Vert_{L_2(\Omega)}^2\right]
    \end{equation}
    where $C$ is a constant depending on $d$, $\bar b$, $a_0$ and norms of coefficients $a_i,b_i,c_i,r$ in respective spaces given in \eqref{databounds}.
\end{corollary}
\emph{Proof.} The uniqueness is proved in Proposition~\ref{some}. The existence of the weak solution is a direct consequence of Theorem \ref{approx}. Indeed, given $f\in L_\infty(D)$, consider the sequence of discrete vectors $[f]_{\D}:=\n Q_{\D}(f)$. Corresponding sequence of interpolations $\n P_{\D}([f]_{\D})$ converge strongly to $f$ in $L_2(D)$, and Theorem \ref{approx} implies the existence of the weak solution $v(x,t; f)\in\overset{\circ}{W}{}_2^{1,1}(D)\cap L_{\infty}(D)$. There is a sequence of multilinear interpolations $\{V_{\D}'\}$ of the solution to the discrete PDE problem, which converge to $v$ weakly in $W_2^{1,1}(D)$, strongly in $L_2(D)$, and pointwise a.e. on $D$. From the discrete maximum estimate \eqref{boundedest} of Theorem~\ref{boundedness} it follows that $\Vert V_{\D}'\Vert_{L_{\infty}(D)}$ is bounded above by the right-hand side of (\ref{boundedest}). Noting that, 
\begin{equation}\label{frelation}
\Vert[f]_{\D}\Vert_{\ell_{\infty}}=\Vert f^{\D}\Vert_{L_\infty(D)}\leq \Vert f\Vert_{L_\infty(D)},
\end{equation}
from \eqref{boundedest}, \eqref{boundv} follows. To prove the energy estimate \eqref{energyv} we use the following two estimates proved in \cite{Abdulla7} ((4.17),(4.18)):
\begin{equation}
\Vert D_xV_{\D}'\Vert_{L_2(D)}^2\leq2^{d+1}T\max\limits_{0\leq k\leq n}\sum\limits_{\n A}h^d\sum\limits_{i=1}^d|v_{\gamma x_i}(k)|^2.\label{pwlinearxok}
\end{equation}
\begin{equation}\label{pwlineartok}
\left\Vert\frac{\partial}{\partial t}V_{\D}'\right\Vert_{L_2(D)}^2\leq2^d \sum\limits_{k=1}^n\tau\sum\limits_{\n A}h^d|v_{\gamma\bar t}(k)|^2\,dx.
\end{equation}
Weak convergence in $W_2^{1,1}(D)$ implies that
\begin{equation}\label{normweakconv}
\Vert D_xv\Vert_{L_2(D)}\leq\liminf_{\D\ra0}\Vert D_xV_{\D}'\Vert_{L_2(D)}, \ \  \Vert v_t\Vert_{L_2(D)}\leq\liminf_{\D\ra0}\Big\Vert \frac{\partial}{\partial t}V_{\D}'\Big\Vert_{L_2(D)}.
\end{equation}
From \eqref{pwlinearxok},\eqref{pwlineartok}, \eqref{normweakconv}, \eqref{phibound}, \eqref{frelation} and \eqref{energyest}, \eqref{energyv} follows.\hfill{$\square$}

Having estimates \eqref{boundedest},\eqref{energyest}, and approximation Theorem~\ref{approx}, the completion of the proofs of Theorems~\ref{optsol} and \ref{funcapprox} coincides with the proofs given in \cite{Abdulla7}. Theorem~\ref{approx} implies that the cost functional $\n J(f)$ is continuous on $\n F^R$ in a weak topology of $L_2(D)$. Therefore, existence of the optimal control is a consequence of the Weierstrass theorem in a weak topology due to weak compactness of the control set $\n F^R$ \cite{Fursikov}. Proof of the  convergence with respect to functional, or claim\eqref{approx1} of Theorem~\ref{funcapprox} is pursued by proving claims (i) and (ii) of the Lemma~\ref{Vasil}. Claim of Theorem~\ref{funcapprox} on the convergence with respect to control is a direct consequence of Theorem~\ref{approx}.\\

\end{document}